\documentclass[3p,times]{elsarticle}
\usepackage{amsmath, amssymb, amsthm, amscd}
\usepackage{amsfonts}
\usepackage{graphicx}
\usepackage{color}
\usepackage{mathrsfs}
\usepackage{framed}
\usepackage[breaklinks,colorlinks]{hyperref}

\makeatletter
\newcommand{\autorefcheckize}[1]{%
  \expandafter\let\csname @@\string#1\endcsname#1%
  \expandafter\DeclareRobustCommand\csname relax\string#1\endcsname[1]{%
    \csname @@\string#1\endcsname{##1}\wrtusdrf{##1}}%
  \expandafter\let\expandafter#1\csname relax\string#1\endcsname
}
\makeatother

\usepackage{setspace}

\newtheorem{theorem}{Theorem}

\newtheorem{lemma}[theorem]{Lemma}

\newcommand{\ra}{\rightarrow}
\newcommand{\p}{\partial}
\newcommand{\f}{\frac}

\renewcommand{\f}{\frac}

\newcommand{\be}{\begin{eqnarray}}
\newcommand{\ee}{\end{eqnarray}}
\newcommand{\bea}{\begin{eqnarray}}
\newcommand{\eea}{\end{eqnarray}}
\newcommand{\bna}{\begin{eqnarray*}}
\newcommand{\ena}{\end{eqnarray*}}

\renewcommand{\le}{\left}
\newcommand{\ri}{\right}

\usepackage{bbding}

\allowdisplaybreaks

\journal{XXX}
\begin{document}

\begin{frontmatter}

\title{Boundary value problem for the mean field equation on a compact Riemann surface \tnoteref{sw}}

\author{Jiayu Li $^{1}$\corref{y}}
\cortext[y]{Corresponding author.}
 \ead{jiayuli@ustc.edu.cn}

\author{Linlin Sun $^2$}
 \ead{ sunll@whu.edu.cn}

\author{Yunyan Yang $^3$}
 \ead{yunyanyang@ruc.edu.cn}

\address{$^1$ School of Mathematical Sciences,
University of Science and Technology of China, Hefei, 230026, China}
\address{$^2$ School of Mathematics and Statistics, Wuhan University, Wuhan 430072, China}
 \address{$^3$ Department of Mathematics,
Renmin University of China, Beijing 100872, China}

\tnotetext[sw]{This research is partly supported by the National Natural Science Foundation of China (Grant No.11721101),
and by the National Key Research and Development Project SQ2020YFA070080.}

\begin{abstract}

Let $(\Sigma,g)$ be a compact Riemann surface with smooth boundary $\p\Sigma$, $\Delta_g$ be the Laplace-Beltrami operator,
and $h$ be a positive smooth function. Using a min-max scheme introduced by Djadli-Malchiodi (2006) and Djadli
(2008), we prove that if $\Sigma$ is non-contractible, then
 for any $\rho\in(8k\pi,8(k+1)\pi)$ with $k\in\mathbb{N}^\ast$, the mean field equation
$$\le\{\begin{array}{lll}
\Delta_g u=\rho\f{he^u}{\int_\Sigma he^udv_g}&{\rm in}&\Sigma\\[1.5ex]
u=0&{\rm on}&\p\Sigma
\end{array}\ri.$$
has a solution. This generalizes earlier existence results of Ding-Jost-Li-Wang (1999) and Chen-Lin (2003) in the Euclidean
domain.

Also we consider the corresponding Neumann boundary value problem. If $h$ is a positive smooth function, then
for any $\rho\in(4k\pi,4(k+1)\pi)$ with $k\in\mathbb{N}^\ast$,
the mean field equation
$$\le\{\begin{array}{lll}
\Delta_g u=\rho\le(\f{he^u}{\int_\Sigma he^udv_g}-\f{1}{|\Sigma|}\ri)&{\rm in}&\Sigma\\[1.5ex]
\p u/\p{\mathbf{v}}=0&{\rm on}&\p\Sigma
\end{array}\ri.$$
has a solution, where $\mathbf{v}$ denotes the unit normal outward vector on $\p\Sigma$. Note that in this case we do not
require the surface to be non-contractible.

\end{abstract}


\end{frontmatter}

\section{Introduction}
As a basic problem of mathematical physics, the mean field equation has aroused the interests of many mathematicians for at least half a century. In addition to the prescribed Gaussian curvature problem \cite{Berger,Chang-Yang1,Chang-Yang2,Chen-Ding,KW}, it also arises in Onsager's vortex model for turbulent Euler flows \cite[Page 256]{M-P}, and in Chern-Simons-Higgs models \cite{C-Y,DJLW4,DJLW3,S-T,Tarantello,YY}.

Let $\Omega$ be a smooth bounded domain in $\mathbb{R}^2$. It was proved by Ding-Jost-Li-Wang \cite{DJLW2} that if the complement of $\Omega$ contains a bounded region, and $h:\overline\Omega\ra\mathbb{R}$ is a positive function, then the mean field equation
\be\label{Dirichlet-domain}
\le\{\begin{array}{lll}
-\Delta_{\mathbb{R}^2} u=\rho\f{he^u}{\int_\Omega he^udx}&{\rm in}&\Omega\\[1.5ex]
u=0&{\rm on}&\p\Omega
\end{array}\ri.
\ee
has a solution for all $\rho\in(8\pi,16\pi)$, where $\Delta_{\mathbb{R}^2}=\p^2/\p x_1^2+\p^2/\p x_2^2$ is the standard Laplacian operator in $\mathbb{R}^2$. The proof is based on a compactness result of Li-Shafrir \cite{Li-Shafrir}, the monotonicity technique used by Struwe \cite{Struwe} in dealing with harmonic maps, and a general min-max theorem \cite[Theorem 2.8]{Willem}.

It was pointed out by Li \cite{Li} that the Leray-Schauder degree for the mean field equation should depend only on
the topology of the domain and $k\in\mathbb{N}$ satisfying $\rho\in(8k\pi,8(k+1)\pi)$. To illustrate this point, he calculated the simplest case $\rho<8\pi$.  Later, by computing the topological degree, Chen-Lin \cite{Chen-Lin} improved Ding-Jost-Li-Wang's result to the following: If $\Omega$ is not simply connected, and $h$ is positive on $\overline\Omega$, then \eqref{Dirichlet-domain} has a solution for all $\rho\in(8k\pi,8(k+1)\pi)$. Also they were able to compute the topological degree for the mean field equation on compact Riemann surface $(\Sigma,g)$ without boundary, namely
\be\label{surface}
\Delta_{g} u=\rho\le(\f{he^u}{\int_\Sigma he^udv_g}-\f{1}{|\Sigma|}\ri)\quad{\rm in}\quad\Sigma,
\ee
where $\Delta_g$ denotes the Laplace-Beltrami operator, and $|\Sigma|$ stands for the area of $\Sigma$ with respect to the metric $g$. Precisely the degree-counting formula for \eqref{surface} is given by
$\binom{k-\chi(\Sigma)}{k}$ for $\rho\in(8k\pi,8(k+1)\pi)$. As a consequence, if the Euler characteristic $\chi(\Sigma)\leq 0$, then \eqref{surface} has a solution.

 Note that solutions of \eqref{surface} are critical points of the functional
 $$J_\rho(u)=\f{1}{2}\int_\Sigma|\nabla_gu|^2dv_g-\rho\log\int_\Sigma he^udv_g+\f{\rho}{|\Sigma|}\int_\Sigma udv_g,\quad u\in W^{1,2}(\Sigma).$$
 A direct method of variation leads to that $J_\rho$ has critical points for $\rho<8\pi$. When $\rho=8\pi$, Ding-Jost-Li-Wang \cite{DJLW} find a critical point of $J_{8\pi}$ under certain conditions on $\Sigma$ and $h$. For the cases $\rho=8\pi$, $h\geq 0$ or $h$ changes sign, we refer the readers to \cite{Sun-Zhu2,Sun-Zhu,Yang-Zhu}.

 In a celebrated paper, Djadli \cite{Djadli} were able to find
 a solution of \eqref{surface} for all $\rho\in(8k\pi,8(k+1)\pi)$ ($k\in\mathbb{N}^\ast$) and arbitrary genus of $\Sigma$, by adapting a min-max scheme introduced by Djadli-Malchiodi \cite{Djadli-Malchiodi}. In particular Chen-Lin's existence result for \eqref{surface} was improved by Djadli to arbitrary possible $\chi(\Sigma)$. Let us summarize the procedure in \cite{Djadli}. Denote the family of formal sums by
 $$\Sigma_k=\le\{\sum_{i=1}^kt_i\delta_{x_i}: t_i\geq 0,\,\sum_{i=1}^kt_i=1,\,x_i\in \Sigma\ri\},$$
 endowed with the weak topology of distributions, say the topology of $(C^1(\Sigma))^\ast$. This is known in literature as the formal set of barycenters of $\Sigma$. The {\it first} and most important step is to construct two continuous map $\Psi$ and $\Phi_\lambda$ between $\Sigma_k$ and sub-levels of $J_\rho$, say
$$\Sigma_k\overset{\Phi_{\lambda}}{\ra}J_{\rho}^{-\left(\rho-8k\pi\right)\ln\lambda}\overset{\Psi}{\ra}\Sigma_k$$
for $\lambda\geq\lambda_L=e^{L/(\rho-8k\pi)}$ and large $L>0$; moreover $\lim_{\lambda\to+\infty}\Psi\circ\Phi_{\lambda}=\mathrm{Id}$, in particular, $\Psi\circ\Phi_{\lambda}$ is homotopic to the identity on $\Sigma_k$ provided $\lambda\geq\lambda_L$. Here $J_\rho^{a}$ stands for a set of all functions $u\in W^{1,2}(\Sigma)$ with $J_\rho(u)\leq a$ for any real number $a$. The {\it second} step is to set suitable min-max value for $J_\rho$, namely
$$\alpha_{\lambda,\rho}=\inf_{\gamma\in\Gamma_\lambda}\sup_{(\sigma,t)\in\widehat{\Sigma}_k} J_\rho(\gamma(\sigma,t)),$$
where $\widehat{\Sigma}_k=\Sigma_k\times[0,1]/(\Sigma_k\times\{0\})$ is a topological cone, and $\Gamma_\lambda$ is a set of paths
$$\Gamma_\lambda=\le\{\gamma\in C^0(\widehat{\Sigma}_k,W^{1,2}(\Sigma)): \gamma(\sigma,1)=\Phi_\lambda(\sigma),\,\forall\sigma\in\Sigma_k\ri\}.$$
The hypothesis $\rho\in(8k\pi,8(k+1)\pi)$ and the fact that $\Sigma_k$ is non-contractible lead to $\alpha_{\lambda,\rho}>-\infty$ for sufficiently large $\lambda$. The {\it third} step is to get critical points of $J_\rho$ for $\rho\in\Lambda$, where $\Lambda$ is a dense subset of $(8k\pi,8(k+1)\pi)$, by using the monotonicity of $\alpha_{\lambda,\rho}/\rho$. The {\it final} step is to find  critical points of $J_\rho$ for any $\rho\in(8k\pi,8(k+1)\pi)$, by using a compactness result of Li-Shafrir \cite{Li-Shafrir} and an improved Trudinger-Moser inequality due to Chen-Li \cite{Chen-Li}.
Note that the last two steps are essentially done by Ding-Jost-Li-Wang \cite{DJLW2}.

This method was extensively used to deal with  the problems of elliptic equations or systems involving exponential growth nonlinearities. For Toda systems, we refer the readers to \cite{Ao,BJMR,Malchiodi-Nd,Malchiodi-Ruiz} and the references therein. Recently Sun-Wang-Yang \cite{Sun} extended Djadli's result to the case of a generalized mean field equation. Marchis-Malchiodi-Martinazzi-Thizy \cite{MMMT} employed it to find critical points of a Trudinger-Moser functional.

\vspace{2ex}

In this paper, we concern the boundary value problems on the mean field equation. From now on, we let $(\Sigma,g)$ be a compact Riemann surface with smooth boundary $\p\Sigma$. Our first aim is to generalize results of Ding-Jost-Li-Wang \cite{DJLW2} and Chen-Lin \cite{Chen-Lin}. Precisely we have the following:

\begin{theorem}\label{thm1}
Let $(\Sigma,g)$ be a compact Riemann surface with smooth boundary $\p\Sigma$, $\Delta_g$ be the Laplace-Beltrami operator, and $h:\overline\Sigma\ra \mathbb{R}$ be a positive smooth function. If $\Sigma$ is not  simply connected, then for any $\rho\in(8k\pi,8(k+1)\pi)$ with $k$ a positive integer,  the Dirichlet problem
\be\label{Dir}\le\{\begin{array}{lll}
\Delta_gu=\rho\f{he^u}{\int_\Sigma he^udv_g}&{\rm in}& \Sigma\\[1.5ex]
u=0&{\rm on}&\p\Sigma
\end{array}\ri.\ee
has a solution.
\end{theorem}

The proof of \autoref{thm1} is based on the min-max theorem \cite[Theorem 2.8]{Willem}, which was also used by Ding-Jost-Li-Wang \cite{DJLW2} and Djadli \cite{Djadli}, compactness analysis, and an improved Trudinger-Moser inequality. All of the three parts are quite different from those of \cite{DJLW2,Djadli}. On the choice of the metric space, we use $\widehat{\Sigma}_{\epsilon,k}$ (see \eqref{S-ek} below) instead of $\widehat{\overline\Sigma}_k$ or $\widehat{\Sigma}_k$; On compactness analysis, we use a reflection method different from that of Chen-Lin \cite{Chen-Lin} to show  blow-up phenomenon can not occur on the boundary $\p\Sigma$; Moreover, we need to prove an improved Trudinger-Moser inequality for functions with boundary value zero, nor is it the original one in \cite{Chen-Li}.

\vspace{2ex}

We also consider the Neumann boundary value problem on the mean field equation. In this regards, our second result is the following:

\begin{theorem}\label{thm2}
Let $(\Sigma,g)$ be a compact Riemann surface with a smooth boundary $\p\Sigma$, $\mathbf{v}$ be the unit normal outward vector on $\p\Sigma$, $\Delta_g$ be the Laplace-Beltrami operator, and $h:\overline\Sigma\ra \mathbb{R}$ be a positive smooth function. If
$\rho\in(4k\pi,4(k+1)\pi)$ with $k$ a positive integer, then the Neumann boundary value problem
\be\label{num}\le\{\begin{array}{lll}
\Delta_gu=\rho\le(\f{he^u}{\int_\Sigma he^udv_g}-\f{1}{|\Sigma|}\ri)&{\rm in}& \Sigma\\[1.5ex]
{\p u}/{\p{\mathbf{v}}}=0&{\rm on}&\p\Sigma
\end{array}\ri.\ee
has a solution.
\end{theorem}

We remark that in \autoref{thm2}, $\Sigma$ need not to be non-contractible. For the proof of \autoref{thm2}, we choose a metric space $\widehat{\mathscr{S}}_k$ (see \eqref{subset} below), which is non-contractible, whether $\Sigma$ is non-contractible or not. Concerning the compactness of solutions to \eqref{num}, if it has a sequence of blow-up solutions, then we show that $\rho=4k\pi$ for $k\in\mathbb{N}^\ast$. Also we derive an improved Trudinger-Moser inequality for functions with integral mean zero, which is important in our analysis.

Before ending this introduction, we mention a recent result of Zhang-Zhou-Zhou \cite{Zhang-Zhou-Zhou}. Using the min-max scheme of
Djadli and Djadli-Malchiodi, they obtained the existence of solutions to the equation
\begin{eqnarray*}
\left\{\begin{array}{lll}
\Delta_{g} u=0
 & {\rm in }&  \Sigma \\[1.5ex]
{\partial u}/{\partial \mathbf{v}}=\rho  \frac{h e^{u}}{\int_{\partial \Sigma} h e^{u} \mathrm{d}s_g}
 &{\rm on }& \partial \Sigma
\end{array}\right.
\end{eqnarray*}
for any $\rho\in(2k\pi,2(k+1)\pi)$, $k \in \mathbb{N}^{*}$, and any positive smooth function $h$. This improved an early result of Guo-Liu \cite{Guo-Liu}.

\vspace{2ex}

In the remaining part of this paper, \autoref{thm1} and \autoref{thm2} are proved by the min-max method in  \autoref{Dirichlet0} and \autoref{Neuma} respectively. Throughout this paper, Sequence and subsequence are  not distinguished.
We often denote various constants by the same $C$ from line to line, even in the same line. Sometimes we write
 $C_k$, $C_{k,\epsilon}$, $C(\epsilon),\dots$, to emphasize the dependence of these constants.

\section{The Dirichlet boundary value problem}\label{Dirichlet0}
In this section, \autoref{thm1} is proved. This will be divided into several subsections. In \autoref{2.1}, we analyze the compactness of solutions to the Dirichlet problem \eqref{Dir}; In \autoref{2.2}, we derive an improved Trudinger-Moser inequality for functions $u\in W_0^{1,2}(\Sigma)$; In  \autoref{2.3}, we construct two continuous maps between sub-levels $J_\rho^{-L}$ with sufficiently large $L$ and the topological space $\Sigma_k$; In  \autoref{2.4}, we construct min-max levels of $J_\rho$, and ensure these min-max levels are finite; In \autoref{2.5}, several uniform estimates on min-max levels of $J_\rho$ are obtained; In \autoref{2.6}, adapting the argument of \cite[Lemma 3.2]{DJLW2}, we prove that $J_\rho$ has a critical point for $\rho$ in a dense subset of $(8k\pi,8(k+1)\pi)$; In \autoref{2.7}, using compactness of solutions to the Dirichlet problem  \eqref{Dir}, we conclude that $J_\rho$ has a critical point for any $\rho\in(8k\pi,8(k+1)\pi)$.

\subsection{Compactness analysis}\label{2.1}
Let $(\rho_n)$ be a sequence of numbers tending to $\rho$, $(h_n)$ be a function sequence converging to
$h$ in $C^1(\overline{\Sigma})$, and $(u_n)$ be a sequence of solutions to
\be\label{Dirichlet}\le\{\begin{array}{lll}
\Delta_g u_n=\rho_n\f{h_ne^{u_n}}{\int_\Sigma h_ne^{u_n}dv_g}&{\rm in}&\Sigma\\[1.5ex]
u_n=0&{\rm on}&\p\Sigma.
\end{array}\ri.\ee
Denote $v_n=u_n-\log\int_\Sigma h_ne^{u_n}dv_g$. Then
$\Delta_{g}v_n=\rho_nh_ne^{v_n}$ and $\int_\Sigma h_ne^{v_n}dv_g=1$.

\begin{lemma}\label{compact}
Assume $\rho$ is a positive number and $h$ is a positive function. Up to a subsequence, there holds one of the following alternatives:\\
$(i)$ $(u_n)$ is bounded in $L^\infty(\overline\Sigma)$;\\
$(ii)$ $(v_n)$ converges to $-\infty$ uniformly in $\overline\Sigma$;\\
$(iii)$ there exists a finite singular set $\mathcal{S}=\{p_1,\cdots,p_m\}\subset {\Sigma}$ such that
for any $1\leq j\leq m$, there is a sequence of points $\{p_{j,n}\}\subset\Sigma$ satisfying
$p_{j,n}\ra p_j$, $u_n(p_{j,n})\ra+\infty$, and $v_n$ converges to $-\infty$ uniformly on any compact subset of
$\overline{\Sigma}\setminus\mathcal{S}$ as $n\ra\infty$. Moreover,
$$\rho_n\int_\Sigma h_ne^{v_n}dv_g\ra 8m\pi. $$
\end{lemma}

\proof Note that $\overline{\Sigma}=\Sigma\cup \p\Sigma$, where $\Sigma$ is an open set including all inner points of
$\overline{\Sigma}$,
 and $\p\Sigma$ is its boundary. The compactness analysis on $(u_n)$ will be divided into two parts.\\

{\bf Part I. Analysis in the interior domain $\Sigma$}

According to an observation in \cite[Section 4.1]{Uniform} (compared with \cite[Theorem 4.17]{Aubin}), we have by the Green representation formula for functions with boundary value zero,
\be\label{W1q}\|u_n\|_{W_0^{1,q}(\Sigma)}\leq C_q,\quad\forall 1<q<2.\ee
We claim that there exists some constant $c_0>0$ such that for all $n\in \mathbb{N}$,
\be\label{lower-b}\int_\Sigma h_ne^{u_n}dv_g\geq c_0.\ee
Suppose not. By Jensen's inequality
$$e^{\f{1}{|\Sigma|}\int_\Sigma u_ndv_g}\leq \f{1}{|\Sigma|}\int_\Sigma e^{u_n}dv_g\ra 0.$$
Thus $\int_\Sigma u_ndv_g\ra -\infty$, which contradicts \eqref{W1q}, and concludes our claim \eqref{lower-b}.

To proceed, we assume $\rho_nh_ne^{v_n}dv_g$ converges to some nonnegative measure $\mu$. If $\mu(x^\ast)<4\pi$ for some $x^\ast\in\Sigma$, then there exist two positive constants $\epsilon_0$ and $r_0$ such that
$$\int_{B_{x^\ast}(r_0)}\rho_nh_ne^{v_n}dv_g\leq 4\pi-\epsilon_0.$$
In view of \eqref{Dirichlet}, by Brezis-Merle's theorem \cite[Theorem 1]{B-M} and elliptic estimates, we have that $(u_n)$ is bounded in $L^{\infty}(B_{x^\ast}(r_0/2))$. This leads to $\mu(x^\ast)=0$. Define a set $\mathcal{S}=\{x\in\Sigma: \mu(x)\geq 4\pi\}$.

If $\mathcal{S}\not=\varnothing$, then we shall show that for any compact set $A\subset\Sigma\setminus\mathcal{S}$, there holds
\be\label{local-uniform}v_n\ra-\infty\,\,{\rm uniforlmy\,\,in}\,\, x\in A.\ee
It suffices to prove  that
\be\label{infty}\int_\Sigma h_ne^{u_n}dv_g\ra+\infty.\ee
Suppose \eqref{infty} does not hold. In view of \eqref{lower-b}, there is a constant $c_1$ such that up to a subsequence
 $$0<c_0\leq\int_\Sigma h_ne^{u_n}dv_g\leq c_1.$$
Choose $x_0\in\mathcal{S}$ and $0<r_0<{\rm dist}(x_0,\p\Sigma)$ satisfying
 $B_{x_0}(r_0)\cap \mathcal{S}=\{x_0\}$. Note that $(u_n)$ is locally uniformly bounded in $\Sigma\setminus\mathcal{S}$. There exists a positive constant $c_2$ depending on $x_0$ and $r_0$ such that $|v_n(x)|\leq c_2$ for all $x\in \p B_{x_0}(r_0)$.
 Let $w_n$ be a solution to
 $$\le\{
 \begin{array}{lll}
 \Delta_gw_n=\rho_nh_ne^{v_n}&{\rm in}& B_{x_0}(r_0)\\[1.5ex]
 w_n=-c_2&{\rm on}& \p B_{x_0}(r_0).
 \end{array}
 \ri.$$
 Then the maximum principle implies that $w_n\leq v_n$ in $B_{x_0}(r_0)$. By the Green formula, $w_n$ converges to $w$ weakly in $W^{1,q}(B_{x_0}(r_0))$ and a.e. in $B_{x_0}(r_0)$. Moreover, $w$ is a solution of
 $$\le\{
 \begin{array}{lll}
 \Delta_gw=\mu&{\rm in}& B_{x_0}(r_0)\\[1.5ex]
 w=-c_2&{\rm on}& \p B_{x_0}(r_0).
 \end{array}
 \ri.$$
 Let $G_{x_0}$ be a distributional solution of
 $$\le\{
 \begin{array}{lll}
 \Delta_gG_{x_0}=4\pi \delta_{x_0}&{\rm in}& B_{x_0}(r_0)\\[1.5ex]
 G_{x_0}=-c_2&{\rm on}& \p B_{x_0}(r_0).
 \end{array}
 \ri.$$
 Clearly $G_{x_0}$ is represented by
 \be\label{Green}G_{x_0}(x)=-2\log{\rm dist}(x,x_0)+A_{x_0}+o(1),\ee
 where $A_{x_0}$ is a constant, and $o(1)\ra 0$ as $x\ra x_0$. Since
 $$\le\{\begin{array}{lll}
 \Delta_g(w-G_{x_0})\geq 0&{\rm in}& B_{x_0}(r_0)\\[1.5ex]
 w-G_{x_0}=0&{\rm on}& \p B_{x_0}(r_0).
 \end{array}
 \ri.$$
 It follows from the maximum principle that
 \be\label{less}w(x)\geq G_{x_0}(x)\,\,{\rm for\,\,all}\,\, x\in B_{x_0}(r_0)\setminus \{x_0\}.\ee
 Combining \eqref{Green}, \eqref{less} and the fact that $w_n\ra w$ a.e. in $B_{x_0}(r_0)$, by the Fatou Lemma, we calculate
 \begin{equation*}
 +\infty=\int_{B_{x_0}(r_0)}e^{G_{x_0}}dv_g\leq\int_{B_{x_0}(r_0)}e^wdv_g\leq \liminf_{n\ra\infty}\int_{B_{x_0}(r_0)}e^{w_n}dv_g\leq \liminf_{n\ra\infty}\int_{B_{x_0}(r_0)}
 e^{v_n}dv_g\leq 1.
 \end{equation*}
 This is impossible and excludes the possibility of \eqref{infty}. Hence we conclude \eqref{local-uniform}.

  We may assume  $\mathcal{S}=\{x_1,\cdots,x_m\}$. Then there would hold $\mu(x_i)=8\pi$ for all $1\leq i\leq m$.
With no loss of generality, it suffices to prove $\mu(x_1)=8\pi$. Choose an isothermal coordinate system
$\phi:U\ra\mathbb{B}_1(0)$ near $x_1$. In such coordinates, the metric $g$ and the Laplace-Beltrami operator $\Delta_g$
are represented by $g=e^{\psi(y)}(dy_1^2+dy_2^2)$ and
 $\Delta_g=-e^{-\psi(y)}\Delta_{\mathbb{R}^2}$ respectively, where $\psi$ is a smooth function with $\psi(0,0)=0$, and $\Delta_{\mathbb{R}^2}={\p^2}/{\p y_1^2}+{\p^2}/{\p y_2^2}$ denotes the standard Laplacian on $\mathbb{R}^2$.
 Set $\widetilde{u}=u\circ \phi^{-1}$ for any function $u:U\ra\mathbb{R}$. Since $(u_n)$ is a sequence of solutions to \eqref{Dirichlet},
 $\widetilde{u}_n=u_n\circ \phi^{-1}$ satisfies
 \be\label{local}-\Delta_{\mathbb{R}^2}\widetilde{u}_n(y)=e^{\psi(y)}\rho_n\widetilde{h}_n(y)e^{\widetilde{v}_n(y)},
 \quad y\in\mathbb{B}_1.\ee
 Multiplying both sides of \eqref{local} by $y\cdot\nabla_{\mathbb{R}^2}\widetilde{u}_n(y)$, we have by integration by parts
 \bea\label{Pohozaev}\nonumber\f{r}{2}\int_{\p\mathbb{B}_r}|\nabla_{\mathbb{R}^2}\widetilde{u}_n|^2d\sigma-r\int_{\p\mathbb{B}_r}
 \langle\nabla_{\mathbb{R}^2}\widetilde{u}_n,\mathbf{\nu} \rangle^2d\sigma&=&r\int_{\p\mathbb{B}_r}
 e^{\psi}\rho_n\widetilde{h}_ne^{\widetilde{v}_n}d\sigma-\int_{\mathbb{B}_r}e^{\widetilde{v}_n}
 \rho_n\langle\nabla_{\mathbb{R}^2}(e^{\psi}\widetilde{h}_n),y \rangle dy\\&&-2\int_{\mathbb{B}_r}
 e^{\psi}\rho_n\widetilde{h}_ne^{\widetilde{v}_n}dy,\eea
 where $\mathbb{B}_r=\{y\in\mathbb{R}^2:y_1^2+y_2^2<r\}$, $\p\mathbb{B}_r=\{y\in\mathbb{R}^2:y_1^2+y_2^2=r\}$, and
 $\nu$ denotes the unit outward vector on $\p\mathbb{B}_r$. In view of \eqref{local-uniform}, $(u_n)$ converges to
 a Green function $G(x,\cdot)$ weakly in $W_0^{1,q}(\Sigma)$, and in $C^2_{\rm loc}(\Sigma\setminus\mathcal{S})$. Locally
 $G(x_1,\cdot)$ satisfies
 $$\Delta_{g,z}G(x_1,z)=\mu(x_1)\delta_{x_1}(z),\quad\forall z\in \phi^{-1}(\mathbb{B}_{1}).$$
 Clearly $\widetilde{G}(y)=G(x_1,\phi^{-1}(y))=-\f{\mu(x_0)}{2\pi}\log|y|+\eta(y)$ for some $\eta\in C^2(\mathbb{B}_1)$.
 Passing to the limit $n\ra\infty$ first, and then $r\ra 0$ in \eqref{Pohozaev}, we obtain
 $$\mu(x_1)=\lim_{r\ra 0}\le(\f{r}{2}\int_{\p\mathbb{B}_r}\langle\nabla_{\mathbb{R}^2}\widetilde{G},\nu\rangle^2d\sigma-\f{r}{4}
 \int_{\p\mathbb{B}_r}|\nabla_{\mathbb{R}^2}\widetilde{G}|^2d\sigma\ri)=\f{(\mu(x_1))^2}{8\pi}.$$
 This immediately leads to $\mu(x_1)=8\pi$.\\

{\bf Part II. Analysis on the boundary $\p\Sigma$}

 Let $x^\ast\in \Sigma$ be fixed. Note that $\rho_nh_ne^{v_n}dv_{g}$ converges to the nonnegative Radon measure $\mu$ on $\overline{\Sigma}$.
 If $\mu(x^\ast)<2\pi$, there exist a neighborhood $V$ of $x^\ast$ and a number $\gamma_0>0$ such that
 \be\label{2pi}\int_V\rho_nh_ne^{v_n}dv_g\leq 2\pi-\gamma_0.\ee
 With no loss of generality, we take an isothermal coordinate system $(V,\phi,\{y_1,y_2\})$ such that
 $\phi(x^\ast)=(0,0)$, and $\phi:V\ra{\mathbb{B}_1^+}\cup\Gamma=\{(y_1,y_2):y_2\geq 0\}$, where
 $\Gamma=\{(y_1,y_2): |y_1|<1, y_2=0\}$. Moreover, in this coordinate system, the metric $g=e^{\psi(y)}(dy_1^2+dy_2^2)$, and the Laplace-Beltrami
 operator $\Delta_g=-e^{-\psi(y)}\Delta_{\mathbb{R}^2}$, where $\psi:\mathbb{B}_1^+\cup\Gamma\ra \mathbb{R}$ is a smooth function
 with $\psi(0,0)=0$. For more details about isothermal coordinates on the boundary, we refer the readers to
 \cite{Yang-Zhou}. Now the local version of \eqref{Dirichlet} reads
 \be\label{semi-equation}\le\{\begin{array}{lll}
 -\Delta_{\mathbb{R}^2}(u_n\circ\phi^{-1})(y)=e^{\psi(y)}\rho_n(h_n\circ\phi^{-1})(y)e^{(v_n\circ\phi^{-1})(y)}
 &{\rm in}&\mathbb{B}_1^+\\[1.5ex]
 u_n\circ\phi^{-1}(y)=0&{\rm on}& \Gamma.
 \end{array}\ri.\ee
 For any function $u:V\ra \mathbb{R}$, we define a function $\widetilde{u}: \mathbb{B}_1^+\cup\Gamma\ra\mathbb{R}$ by
 \be\label{u-tilde}\widetilde{u}(y_1,y_2)=\le\{
 \begin{array}{lll}
 u\circ\phi^{-1}(y_1,y_2)&{\rm if}& y_2\geq 0\\[1.5ex]
 -u\circ\phi^{-1}(y_1,-y_2)&{\rm if}& y_2< 0.
 \end{array}\ri.\ee
 One can easily check that $\widetilde{u}_n$ is a distributional solution of
 \be\label{un}-\Delta_{\mathbb{R}^2}\widetilde{u}_n(y)=\widetilde{f}_n(y),\quad y\in\mathbb{B}_1,\ee
 where $\widetilde{f}_n$ is defined as in \eqref{u-tilde} and for $y\in \mathbb{B}_1^+\cup\Gamma$,
 $$f_n\circ\phi^{-1}(y)=e^{\psi(y)}\rho_n(h_n\circ\phi^{-1})(y)e^{(v_n\circ\phi^{-1})(y)}.$$
 In view of \eqref{2pi} and the fact that $\psi(0,0)=0$, there would exist a number $0<r_0<1$ such that
 $$\int_{\mathbb{B}_{r_0}}|\widetilde{f}_n(y)|dy\leq 4\pi-{\gamma_0}.$$
 Let $w_n$ be a solution of
 $$\le\{\begin{array}{lll}
 -\Delta_{\mathbb{R}^2}w_n=\widetilde{f}_n&{\rm in}&\mathbb{B}_{r_0}\\[1.5ex]
 w_n=0&{\rm on}&\p\mathbb{B}_{r_0}.
 \end{array}\ri.$$
 By Brezis-Merle's theorem \cite[Theorem 1]{B-M}, there exists some constant $C$ depending only on $\gamma_0$ and $r_0$ such that
 $$\int_{\mathbb{B}_{r_0}}\exp\le(\f{(4\pi-\gamma_0/2)|w_n|}{\|\widetilde{f}_n\|_{L^1(\mathbb{B}_{r_0})}}\ri)dy\leq C.$$
 Hence there exists some $q_0>1$ such that
 \be\label{Lq0}\|e^{|w_n|}\|_{L^{q_0}(\mathbb{B}_{r_0})}\leq C.\ee
 Let
 $\eta_n=\widetilde{u}_n-w_n$. Then $\eta_n$ satisfies
 \be\label{harmonic}\le\{\begin{array}{lll}
 -\Delta_{\mathbb{R}^2}\eta_n=0&{\rm in}&\mathbb{B}_{r_0}\\[1.5ex]
 \eta_n=\widetilde{u}_n&{\rm on}&\p\mathbb{B}_{r_0}.
 \end{array}\ri.\ee
 Noticing \eqref{W1q} and \eqref{Lq0}, we have by applying elliptic estimates to \eqref{harmonic} that
 \be\label{etan}\|\eta_n\|_{L^\infty(\mathbb{B}_{r_0/2})}\leq C.\ee
 Combining \eqref{lower-b}, \eqref{Lq0} and \eqref{etan}, we conclude
 $\|\widetilde{f}_n\|_{L^{q_0}(\mathbb{B}_{r_0/2})}\leq C$. Applying elliptic estimates to \eqref{un}, we obtain that 
 $\|\widetilde{u}_n\|_{L^\infty(\mathbb{B}_{r_0/4})}\leq C$, which implies
 $\|{u}_n\|_{L^\infty(\phi^{-1}(\mathbb{B}^+_{r_0/4}))}\leq C$. In conclusion, we have that if $\mu(x^\ast)<2\pi$, then
 $(u_n)$ is uniformly bounded near $x^\ast$. This also leads to $\mu(x^\ast)=0$.

 If  $\mu(x^\ast)\geq 2\pi$, in the same coordinate system $(V,\phi,\{y_1,y_2\})$ as above, $\widetilde{f}_n(y)dy$ converges to a Radon measure $\widetilde{\mu}$ with
 $\widetilde{\mu}(0,0)\geq 4\pi$. Obviously there exists some $r_1>0$ such that for any $x\in \mathbb{B}_{r_1}\setminus\{(0,0)\}$,
 $\widetilde{\mu}(x)=0$. Using the same argument as the proof  of \eqref{local-uniform}, we conclude that
 for any compact set $A\subset \mathbb{B}_{r_1}\setminus\{(0,0)\}$, $\widetilde{v}_n$ converges to $-\infty$ uniformly in $A$. This leads to
 $\widetilde{f}_n(y)dy$ converges to the Dirac measure $\widetilde{\mu}(0,0)\delta_{(0,0)}(y)$. Recalling \eqref{W1q}, we have
 $\widetilde{u}_n$ converges to $\widetilde{G}_0$ weakly in $W^{1,q}(\mathbb{B}_{r_1})$ and a.e. in $\mathbb{B}_{r_1}$, where $\widetilde{G}_0$
 satisfies
 $$-\Delta_{\mathbb{R}^2}\widetilde{G}_0(y)=\widetilde{\mu}(0,0)\delta_{(0,0)}(y),\quad y\in\mathbb{B}_{r_1}.$$
 Clearly $\widetilde{G}_0$ is represented by
 \be\label{G-rep}\widetilde{G}_0(y)=-\f{\widetilde{\mu}(0,0)}{2\pi}\log|y|+A_0+O(|y|)\ee
 as $y\ra 0$, where $A_0$ is a constant. Noting that $\widetilde{v}_n$ converges to $-\infty$
 locally uniformly in $\mathbb{B}_{r_1}\setminus\{(0,0)\}$, we have by applying elliptic estimates to \eqref{un} that
 \be\label{C1}\widetilde{u}_n\ra \widetilde{G}_0\quad {\rm in}\quad C^1_{\rm loc}(\mathbb{B}_{r_1}\setminus\{(0,0)\}).\ee
 By \eqref{semi-equation}, $\widetilde{u}_n(y_1,0)=0$ for all $|y_1|<1$, which together with \eqref{C1} leads to
 $\widetilde{G}_0(y_1,0)=0$ for all $0<|y_1|<r_1$. This contradicts \eqref{G-rep}. Therefore
 $$\{x\in\p\Sigma:\mu(x)\geq 2\pi\}=\varnothing.$$

 Combining Parts I and II, we conclude the lemma. $\hfill\Box$

 \subsection{An improved Trudinger-Moser inequality}\label{2.2}

 In this subsection, we shall derive an improved Trudinger-Moser inequality, which is analog of that of
 Chen-Li \cite{Chen-Li}. It is known (see for example Jiang \cite{Jiang}) that
 \be\label{T-M}\log\int_\Sigma e^udv_g\leq \f{1}{16\pi}\int_\Sigma |\nabla_gu|^2dv_g+C,\quad\forall u\in
 W_0^{1,2}(\Sigma).\ee

 \begin{lemma}\label{improved}
 Let $b_0>0$ and $\gamma_0>0$ be two constants, $\Omega_1,\cdots,\Omega_k$ be $k$ domains of $\overline\Sigma$ satisfying
 ${\rm dist}(\Omega_i,\Omega_j)\geq b_0$ for all $1\leq i<j\leq k$. Then for any $\epsilon>0$, there exists some constant $C$
  depending only on $b_0,\gamma_0,k,\epsilon$, such that
 \be\label{imp}\log\int_\Sigma e^udv_g\leq \f{1}{16k\pi-\epsilon}\int_\Sigma|\nabla_gu|^2dv_g+C\ee
 for all $u\in W_0^{1,2}(\Sigma)$ with
 \be\label{comp}\int_{\Omega_i}e^{u}dv_g\geq\gamma_0\int_\Sigma e^udv_g,\,\, i=1,\cdots,k.\ee
 \end{lemma}

 \proof  We modify an argument of Chen-Li \cite{Chen-Li}. Take smooth functions $\phi_1,\cdots,\phi_k$
 defined on $\overline\Sigma$ satisfying
 \be\label{phi-1}{\rm supp}\phi_i\cap {\rm supp}\phi_j=\varnothing,\quad\forall 1\leq i<j\leq k,\ee
 \be\label{phi-2}\phi_i\equiv 1 \,\,{\rm on}\,\,\Omega_i; \,\,0\leq\phi_i\leq 1 \,\,{\rm on}\,\, \overline\Sigma,\,\,\forall 1\leq i\leq k, \ee
 and for some positive constant $b_1$ depending only on $b_0$ and $g$,
 \be\label{phi-3}|\nabla_g\phi_i|\leq b_1,\quad\forall 1\leq i\leq k.\ee
 For any $u\in W_0^{1,2}(\Sigma)$ satisfying \eqref{comp}, we have $\phi_i u\in W_0^{1,2}(\Omega)$ for all $1\leq i\leq k$, and thus
 \eqref{T-M} implies
  \bna
  \int_\Sigma e^udv_g&\leq&\f{1}{\gamma_0}\int_{\Omega_i}e^{u}dv_g\\
  &\leq&\f{1}{\gamma_0}\int_{\Sigma}e^{\phi_iu}dv_g\\
  &\leq&\f{1}{\gamma_0}\exp\le(\f{1}{16\pi}\|\nabla_g(\phi_iu)\|_{L^2(\Sigma)}^2+C\ri).
  \ena
  Recall an elementary inequality: if $a\leq a_i$ for nonnegative numbers $a$ and $a_i$, $i=1,\cdots,k$, then
  $a\leq (a_1\cdots a_k)^{1/k}$. In view of \eqref{phi-1}, \eqref{phi-2} and \eqref{phi-3}, we have 
  \bea\nonumber
  \int_\Sigma e^udv_g&\leq&\f{1}{\gamma_0}\le(\prod_{i=1}^k\exp\le(\f{1}{16\pi}\|\nabla_g(\phi_iu)\|_{L^2(\Sigma)}^2+C\ri)\ri)^{1/k}\\
  \nonumber&=&\f{e^C}{\gamma_0}\exp\le(\f{1}{16k\pi}\sum_{i=1}^k\|\nabla_g(\phi_iu)\|_{L^2(\Sigma)}^2\ri)\\\nonumber
  &=&\f{e^C}{\gamma_0}\exp\le(\f{1}{16k\pi}\le\|\nabla_g\le(u\sum_{i=1}^k\phi_i\ri)\ri\|_{L^2(\Sigma)}^2\ri)\\\label{uL2}
  &\leq& C\exp\le(\f{1}{16k\pi}(1+\epsilon_1)\|\nabla_gu\|_{L^2(\Sigma)}^2+C(\epsilon_1)\|u\|_{L^2(\Sigma)}^2\ri).
  \eea

  Let $0<\lambda_1\leq \lambda_2\leq \cdots\leq \lambda_\ell\leq\lambda_{\ell+1}\leq\cdots$ be all eigenvalues of the Laplace-Beltrami operator with respect to the Dirichlet boundary condition with $\lambda_i\ra+\infty$ as $i\ra\infty$, $\{e_i\}_{i=1}^\infty$ be the corresponding unit normal eigenfunctions, i.e., $\Delta_g e_i=\lambda_ie_i$, $\int_\Sigma e_ie_jdv_g=\delta_{ij}$ for $i,j=1,2,\cdots$. It is known that $W_0^{1,2}(\Sigma)=E_\ell\oplus E_\ell^\perp$, where $E_\ell={\rm span}\{e_1,\cdots,e_\ell\}$ and $E_\ell^\perp=\{e_{\ell+1},e_{\ell+2},
  \cdots\}$.
  Let $u\in W_0^{1,2}(\Sigma)$ be as above. Write $u=v+w$ with $v\in E_\ell$ and $w\in E_\ell^\perp$. Thus the Poincar\'e inequality implies
  $$\|v\|_{C^0(\overline\Sigma)}\leq \sum_{i=1}^\ell \|e_i\|_{C^0(\overline\Sigma)}\int_\Sigma |u||e_i|dv_g\leq C_\ell\|\nabla_gu\|_{L^2(\Sigma)}.$$
  While by the definition of the $(\ell+1)$-th eigenvalue,
  $$\int_\Sigma w^2dv_g\leq \f{1}{\lambda_{\ell+1}}\int_\Sigma|\nabla_gw|^2dv_g.$$
  Having the above two estimates and applying \eqref{uL2} to $w$, we have
  \bna
  \int_\Sigma e^udv_g&\leq& e^{C_\ell\|\nabla_gu\|_{L^2(\Sigma)}}\int_\Sigma e^wdv_g\\
  &\leq& Ce^{C_\ell\|\nabla_gu\|_{L^2(\Sigma)}}\exp\le(\f{1}{16k\pi}(1+\epsilon_1)\|\nabla_gw\|_{L^2(\Sigma)}^2+
  \f{C(\epsilon_1)}{\lambda_{\ell+1}}\|\nabla_gw\|_{L^2(\Sigma)}^2\ri)\\
  &\leq&Ce^{C_\ell\|\nabla_gu\|_{L^2(\Sigma)}}\exp\le(\f{1}{16k\pi}\le(1+\epsilon_1\f{C(\epsilon_1)}{\lambda_{\ell+1}}\ri)
  \|\nabla_gu\|_{L^2(\Sigma)}^2\ri).
  \ena
  This together with the Young inequality gives
  \be\label{sec}\log\int_\Sigma e^udv_g\leq \f{1}{16k\pi}\le(1+\epsilon_1\f{C(\epsilon_1)}{\lambda_{\ell+1}}+\epsilon_1\ri)
  \int_\Sigma |\nabla_gu|^2dv_g+C_{\ell,k,\epsilon_1}.
 \ee
 Let $\epsilon>0$ be any given number. Choosing $\epsilon_1=\epsilon/(32k\pi-2\epsilon)$, and then taking a sufficiently large $\ell$ such that  ${C(\epsilon_1)}/{\lambda_{\ell+1}}\leq 1$ in \eqref{sec}, we immediately have
 $$\log\int_\Sigma e^udv_g\leq \f{1}{16k\pi-\epsilon}\int_\Sigma |\nabla_gu|^2dv_g+C,$$
 where $C$ is a constant depending only on $b_0$, $\gamma_0$, $k$, and $\epsilon$. This is exactly \eqref{imp}. $\hfill\Box$\\

 For any $\rho>0$, we define a functional $J_\rho: W_0^{1,2}(\Sigma)\ra\mathbb{R}$ by
 \be\label{functional-0}J_\rho(u)=\f{1}{2}\int_\Sigma|\nabla_gu|^2dv_g-\rho\log\int_\Sigma he^udv_g.\ee
 Clearly $J_\rho\in C^2(W_0^{1,2}(\Sigma),\mathbb{R})$. To find solutions of the mean field equation \eqref{Dir}, it suffices to find critical points of $J_\rho$. For any $a\in \mathbb{R}$,
 the sub-level of $J_\rho$ is written as $J_\rho^{a}=\{u\in W_0^{1,2}(\Sigma): J_\rho(u)\leq a\}$.

 Let $\Sigma_k$ be the formal set of barycenters of $\Sigma$ (of order $k$), which reads as
\be\label{Sig-k}\Sigma_k=\le\{\sum_{i=1}^kt_i\delta_{x_i}:t_i\geq 0,x_i\in\Sigma,\sum_{i=1}^kt_i=1\ri\}.\ee
It is endowed with the weak topology of distributions. In computation, we use on $\Sigma_k$ the metric given by $(C^1(\overline\Sigma))^\ast$ inducing the same topology. Similarly we may define
\be\label{O-S-k}\overline{\Sigma}_k=\le\{\sum_{i=1}^kt_i\delta_{x_i}:t_i\geq 0,x_i\in\overline{\Sigma},\sum_{i=1}^kt_i=1\ri\}.\ee

\subsection{Continuous maps between sub-levels of \texorpdfstring{$J_{\rho}$}{the functional} and \texorpdfstring{$\Sigma_k$}{baycenter}}\label{2.3}

Let $J_\rho$, $\Sigma_k$ and $\overline\Sigma_k$ be defined as in \eqref{functional-0}, \eqref{Sig-k}
and \eqref{O-S-k} respectively. In this subsection, we shall construct continuous maps between sub-levels of $J_\rho$ and $\Sigma_k$ (or $\overline\Sigma_k$).

\begin{lemma}\label{noncontractible}
For any $k\geq 1$, both $\Sigma_k$ and $\overline{\Sigma}_k$ are non-contractible.
\end{lemma}
\proof Since $\Sigma_k$ is homotopic to $\overline{\Sigma}_k$, we only need to prove $\Sigma_k$ is non-contractible.
Let $\chi(\Sigma)$ be the Euler characteristic of $\Sigma$. By \cite[Corollary 1.4 (a)]{Kallel},
\be\label{characteristic}\chi(\Sigma_k)=1-\f{1}{k!}(1-\chi(\Sigma))\cdots(k-\chi(\Sigma)).\ee
Denote the genus $\Sigma$ by  $\mathfrak{g}$ and the connected components of $\partial\Sigma$ by $m$. Notice that $\Sigma$ is  simply connected if and only if the Euler characteristic of $\Sigma$ equal to 1. By the assumption, it then follows  that
\be\label{charac}\chi(\Sigma)=2-2\mathfrak{g}-m\leq 0.\ee
Inserting \eqref{charac} into \eqref{characteristic}, we have
\be\label{char-1}\chi(\Sigma_k)\leq 0.\ee
On the other hand, there holds
\be\label{ch}\chi(\Sigma_k)=\sum_{j=0}^{3k-1}(-1)^j{\rm dim} H_j(\Sigma_k,\mathbb{Z}).\ee
Suppose $\Sigma_k$ is contractible. Then ${\rm dim} H_j(\Sigma_k,\mathbb{Z})=0$ for all $j\geq 1$. While
${\rm dim} H_0(\Sigma_k,\mathbb{Z})=1$, since $\Sigma_k$ is connected. Hence \eqref{ch} gives
$\chi(\Sigma_k)=1$, contradicting \eqref{char-1}. Hence $\Sigma_k$ is non-contractible. $\hfill\Box$

\begin{lemma}\label{Psi}
Let $\rho\in(8k\pi,8(k+1)\pi)$. Then for any sufficiently large $L>0$, there exists a continuous retraction
$$\Psi:J_\rho^{-L}=\le\{u\in W_0^{1,2}(\Sigma): J_\rho(u)\leq -L\ri\}\ra \overline{\Sigma}_k.$$
Moreover, if $(u_n)\subset W_0^{1,2}(\Sigma)$ satisfies
$\f{e^{u_n}}{\int_\Sigma e^{u_n}dv_g}dv_g\ra\sigma\in\overline\Sigma_k$, then $\Psi(u_n)\ra \sigma\in\overline\Sigma_k$.
\end{lemma}

\proof
 By \cite[Proposition 4.1]{Djadli}, for any $\epsilon>0$, there exists $L_0>0$ such that for all $u\in J_\rho^{-L_0}$,
$$\f{e^u}{\int_\Sigma e^udv_g}dv_g\in\le\{\sigma\in\mathcal{D}(\overline\Sigma):{\mathbf d}(\sigma,\overline\Sigma_k)<\epsilon\ri\},$$
where $\mathcal{D}(\overline\Sigma)$ denotes the set of all distributions on $\overline\Sigma$.
If $\epsilon_0>0$ is sufficiently small, then there exists a continuous retraction
\be\label{proj}\psi_k:\le\{\sigma\in\mathcal{D}(\overline\Sigma): {\bf d}(\sigma,\overline{\Sigma}_k)<\epsilon_0\ri\}\ra \overline{\Sigma}_k.\ee
For sufficiently large $L>0$, we set
$$\Psi(u)=\psi_k\le(\f{e^u}{\int_\Sigma e^udv_g}dv_g\ri),\quad \forall u\in J_\rho^{-L}.$$
As a consequence, we have a continuous map $\Psi:J_\rho^{-L}\ra\overline\Sigma_k$.

Moreover, if $(u_n)\subset W_0^{1,2}(\Sigma)$ satisfies
$\f{e^{u_n}}{\int_\Sigma e^{u_n}dv_g}dv_g\ra\sigma\in\overline\Sigma_k$, then as $n\ra\infty$,
$$\Psi(u_n)=\psi_k\le(\f{e^{u_n}}{\int_\Sigma e^{u_n}dv_g}dv_g\ri)\ra \psi_k(\sigma)=\sigma,$$
as desired. $\hfill\Box$\\

Let $\sigma=\sum_{i=1}^kt_i\delta_{x_i}\in\Sigma_k$ be fixed. Take a smooth increasing function $\eta:\mathbb{R}\ra\mathbb{R}$
satisfying  $\eta(t)=t$ for $t\leq 1$, and $\eta(t)=2$ for $t\geq 2$. Set $\eta_r(t)=r\eta(t/r)$ for $r>0$.
For $\lambda>0$ and $x\in\overline\Sigma$, we define
\be\label{ppp}\widetilde{\phi}_{\lambda,\sigma}(x)=\log\le(\sum_{i=1}^kt_i\f{8\lambda^2}{(1+\lambda^2\eta_r^2
({\rm dist}(x,x_i)))^2}\ri)\ee
and
\be\label{p-1}\phi_{\lambda,\sigma}(x)=\widetilde{\phi}_{\lambda,\sigma}(x)-\log \frac{8 \lambda^2}{(1+4 \lambda^{2}r^{2} )^2 }.\ee

\begin{lemma}\label{P}
Let $\rho\in(8k\pi,8(k+1)\pi)$ and $\epsilon>0$. If $\lambda>0$ is chosen sufficiently large, and $r>0$ is chosen sufficiently small, then for any $\sigma\in\Sigma_k$ with ${\rm dist}({\rm supp}\sigma,\p\Sigma)\geq \epsilon$, there hold
\be\label{J-less}J_\rho(\phi_{\lambda,\sigma})\leq \le(8k\pi-\rho\ri)\log\lambda\ee
and
\be\label{tend}\f{e^{\phi_{\lambda,\sigma}}}{\int_\Sigma e^{\phi_{\lambda,\sigma}}dv_g}dv_g\ra \sigma\quad{\rm as}\quad\lambda\ra+\infty.\ee
\end{lemma}
\proof Given $\sigma\in \Sigma_k$. With no loss of generality, we assume ${\rm supp}\sigma=\{x_1,\cdots,x_k\}\subset\Sigma$. Let $\widetilde{\phi}_{\lambda,\sigma}$ and $\phi_{\lambda,\sigma}$ be defined as in \eqref{ppp}
and \eqref{p-1} respectively, where $\lambda>0$ and $0<r<\epsilon/4$.
Write $r_i=r_i(x)={\rm dist}(x,x_i)$ for $x\in\overline\Sigma$. A simple observation gives
\be\label{wide}
\widetilde{\phi}_{\lambda, \sigma}(x)=\left\{\begin{array}{lll}
\log \frac{8 \lambda^2}{(1+4 \lambda^{2}r^{2} )^2 } &{\rm for}& x \in \Sigma\setminus \cup_{i=1}^{k} B_{2 r}(x_i)
\\[1.5ex]
\log \left(\frac{8 \lambda^2 t_{i}}{ \left(1+\lambda^{2} \eta_r^{2}\left(r_i\right) \right)^2 }+\frac{8 \lambda^2\left(1-t_{i}\right)}{ (1+4\lambda^{2} r^{2})^2}\right) &{\rm for}& x \in B_{2 r}(x_i).
\end{array}\right.
\ee
As a consequence, $\phi_{\lambda,\sigma}\in W_0^{1,2}(\Sigma)$.
 For $x\in B_{2r}(x_i)$, a straightforward calculation shows
$$\nabla_g\widetilde\phi_{\lambda,\sigma}(x)=\f{\f{8\lambda^2t_i}{(1+\lambda^2\eta_r^2(r_i))^2}}
{\f{8\lambda^2t_i}{(1+\lambda^2\eta_r^2(r_i))^2}+\f{8\lambda^2(1-t_i)}{(1+4\lambda^2r^2)^2}}
\f{4\lambda^2\eta_r(r_i)\eta_r^\prime(r_i)\nabla_gr_i}{1+\lambda^2\eta_r^2(r_i)},$$
and thus
$$|\nabla_g\widetilde{\phi}_{\lambda,\sigma}(x)|\leq \f{4\lambda^2\eta_r(r_i)\eta_r^\prime(r_i)}{1+\lambda^2\eta_r^2(r_i)}.$$
In view of \eqref{wide}, there holds $\nabla_g\widetilde{\phi}_{\lambda,\sigma}(x)=0$ for $x\in \Sigma\setminus \cup_{i=1}^{k} B_{2 r}(x_i)$. Hence
\bea
\nonumber \int_\Sigma |\nabla_g\widetilde{\phi}_{\lambda,\sigma}|^2dv_g&=&\int_{\cup_{i=1}^kB_{2r}(x_i)}|\nabla_g\widetilde{\phi}_{\lambda,\sigma}|^2dv_g\\
\nonumber&\leq&\sum_{i=1}^k\int_{B_{2r(x_i)}}
\le(\f{4\lambda^2\eta_r(r_i)\eta_r^\prime(r_i)}{1+\lambda^2\eta_r^2(r_i)}\ri)^2dv_g\\
\nonumber&=&\sum_{i=1}^k16\pi(1+O(r^2))\le(\log(1+\lambda^2r^2)+\f{1}{1+\lambda^2r^2}-1\ri)+O(1)\\
&\leq& 16k\pi(1+O(r^2))\log \lambda^2+C\label{nabl}
\eea
for some constant $C$ independent of $r$ and $\lambda$. Moreover, for any $s$ with $0<s<\min\{r,\min_{1\leq i< j\leq k}{\rm dist}(x_i,x_j)\}$, there holds
\bna
\int_{\cup_{i=1}^kB_{2r}(x_i)}e^{\widetilde{\phi}_{\lambda,\sigma}}dv_g&=&\int_{\cup_{i=1}^kB_{s}(x_i)}e^{\widetilde{\phi}_{\lambda,\sigma}}dv_g
+O\le(\f{1}{\lambda^2s^2}\ri)\\
&=&\sum_{i=1}^k\int_{B_{s}(x_i)}
\f{8\lambda^2t_i}{(1+\lambda^2r_i^2)^2}dv_g+O\le(\f{1}{\lambda^2s^2}\ri)\\
&=&8\pi(1+O(s^2))+O\le(\f{1}{\lambda^2s^2}\ri),
\ena
and
\bna
\int_{\Sigma\setminus\cup_{i=1}^kB_{2r}(x_i)}e^{\widetilde{\phi}_{\lambda,\sigma}}dv_g=O\le(\f{1}{\lambda^2r^4}\ri).
\ena
It follows that
\be\label{sense}\int_\Sigma e^{\widetilde\phi_{\lambda,\sigma}}dv_g=8\pi(1+O(s^2))+O\le(\f{1}{\lambda^2s^2}\ri)+O\le(\f{1}{\lambda^2r^4}\ri).\ee
Passing to the limit $\lambda\ra+\infty$ first, and then $s\ra 0+$, we have
\be\label{limit}\lim_{\lambda\ra+\infty}\int_\Sigma e^{\widetilde\phi_{\lambda,\sigma}}dv_g=8\pi.\ee
Hence
\bea\nonumber\int_\Sigma e^{\phi_{\lambda,\sigma}}dv_g&=&\int_\Sigma e^{\widetilde\phi_{\lambda,\sigma}}dv_g \f{(1+4\lambda^2r^2)^2}{8\lambda^2}\\
&=&(8\pi+o_\lambda(1))\lambda^2r^4.\label{phiphi-1}\eea
Combining \eqref{nabl} and \eqref{phiphi-1}, we obtain
\bna
J_\rho(\phi_{\lambda,\sigma})&=&\f{1}{2}\int_\Sigma|\nabla_g\phi_{\lambda,\sigma}|^2dv_g-\rho\log\int_\Sigma he^{\phi_{\lambda,\sigma}}dv_g\\
&\leq&(16k\pi-2\rho+O(r^2))\log\lambda+C_r.
\ena
Since $\rho>8k\pi$, choosing $r>0$ sufficiently small and $\lambda>0$ sufficiently large, we conclude \eqref{J-less}.

Finally we prove \eqref{tend}. Let $\sigma=\sum_{i=1}^kt_i\delta_{x_i}\in\Sigma_k$ be as above. For any $\varphi\in C^1(\overline\Sigma)$, similar to \eqref{sense}, we calculate
$$
\int_\Sigma\varphi e^{\widetilde{\phi}_{\lambda,\sigma}}dv_g
=8\pi\sum_{i=1}^kt_i\varphi(x_i)+O(s^2)+O\le(\f{1}{\lambda^2s^2}\ri)+O\le(\f{1}{\lambda^2r^4}\ri).
$$
Letting $\lambda\ra+\infty$ first, and then $s\ra 0+$, we get
$$\lim_{\lambda\ra+\infty}\int_\Sigma\varphi e^{\widetilde{\phi}_{\lambda,\sigma}}dv_g=8\pi\sum_{i=1}^kt_i\varphi(x_i).$$
This together with \eqref{limit} implies \eqref{tend}.
 $\hfill\Box$

 \begin{lemma}\label{PHI}
 Let $\Psi$ and $L>0$ be as in Lemma \ref{Psi}. If $\lambda>0$ is chosen sufficiently large, then there exists a continuous map $\Phi_\lambda: \Sigma_{k}\ra J_\rho^{-L}$ such that $\Psi\circ \Phi_\lambda$ is homotopic to the identity map ${\rm Id}:\Sigma_k\ra\Sigma_k$.
 \end{lemma}

 \proof Let $\phi_{\lambda,\sigma}$ be constructed as in Lemma \ref{P}. For any $\sigma\in\Sigma_k$, we define
 $\Phi_\lambda(\sigma)=\phi_{\lambda,\sigma}$. Clearly the map $\Phi_\lambda:\Sigma_k\ra W_0^{1,2}(\Sigma)$ is continuous.
 By \eqref{J-less}, if $\lambda\geq e^{L/(\rho-8k\pi)}$, then $J_\rho(\phi_{\lambda,\sigma})\leq -L$. Thus $\Phi_\lambda(\sigma)\in J_\rho^{-L}$. By Lemma \ref{Psi} and \eqref{tend}, there holds
 \bna
 \Psi\circ\Phi_\lambda(\sigma)&=&\Psi(\phi_{\lambda,\sigma})\\
 &=&\psi_k\le(\f{e^{\phi_{\lambda,\sigma}}}{\int_\Sigma e^{\phi_{\lambda,\sigma}}dv_g}dv_g\ri)\\
 &\ra&\sigma
 \ena
 as $\lambda\ra+\infty$. Hence $\Psi\circ \Phi_\lambda$ is homotopic to ${\rm Id}:\Sigma_k\ra\Sigma_k$. $\hfill\Box$

\subsection{Min-max values}\label{2.4}
In this subsection, we shall construct suitable min-max value of $J_\rho$ for $\rho\in(8k\pi,8(k+1)\pi)$, $k\in\mathbb{N}^\ast$.
Recalling that $\Sigma$ is non-contractible, we can take a sufficiently small $\epsilon>0$ such that
$\Sigma_\epsilon=\{x\in \Sigma: {\rm dist}(x,\p\Sigma)\geq \epsilon\}$ is non-contractible.
Let $\Sigma_{\epsilon,k}=\{\sigma\in\Sigma_k: {\rm dist}({\rm supp}\sigma,\p\Sigma)\geq \epsilon\}$. According to Lemma \ref{noncontractible},
we see that $\Sigma_{\epsilon,k}$ is also non-contractible. Let
\be\label{S-ek}\widehat{\Sigma}_{\epsilon,k}=\Sigma_{\epsilon,k}\times [0,1]/(\Sigma_{\epsilon,k}\times\{0\})\ee
 be the topological cone over
$\Sigma_{\epsilon,k}$. A path set associated to
the metric space
$\widehat{\Sigma}_{\epsilon,k}$ is defined by
\be\label{path}\Gamma_\lambda=\le\{\gamma\in C^0(\widehat{\Sigma}_{\epsilon,k},W_0^{1,2}(\Sigma)):\gamma|_{\Sigma_{\epsilon,k}\times\{1\}}\in\Gamma_{\lambda,0}\ri\},\ee
where $\Gamma_{\lambda,0}$ is given by
$$\Gamma_{\lambda,0}=\le\{\gamma\in C^0(\Sigma_{\epsilon,k}\times\{1\},W_0^{1,2}(\Sigma)):\gamma(\sigma,1)=\Phi_\lambda(\sigma),\forall \sigma\in\Sigma_{\epsilon,k}\ri\}.$$
If we write a path $\overline{\gamma}:\widehat{\Sigma}_{\epsilon,k}\ra W_0^{1,2}(\Sigma)$ by $\overline{\gamma}(\sigma,t)=t\phi_{\lambda,\sigma}$,
then $\overline{\gamma}\in \Gamma_\lambda$, and thus $\Gamma_\lambda\not=\varnothing$.

For any real numbers $\lambda$ and $\rho$, we set
\be\label{min-max}\alpha_{\lambda,\rho}=\inf_{\gamma\in\Gamma_\lambda}\sup_{(\sigma,t)\in\widehat{\Sigma}_{\epsilon,k}}J_\rho(\gamma(\sigma,t))\ee
and
\be\label{betarho}\beta_{\lambda,\rho}=\sup_{\gamma\in\Gamma_{\lambda,0}}\sup_{(\sigma,t)\in{\Sigma}_{\epsilon,k}\times\{1\}}
J_\rho(\gamma(\sigma,t)).\ee
\begin{lemma}\label{mini} Let $\rho\in(8k\pi,8(k+1)\pi)$, and $\epsilon>0$ be as above.
If $\lambda$ is chosen sufficiently large, and $r$ is chosen sufficiently small, then
$-\infty<\beta_{\lambda,\rho}<\alpha_{\lambda,\rho}<+\infty$.
\end{lemma}
\proof If $L>0$ is large enough, $\Psi: J_\rho^{-L}\ra\overline\Sigma_k$ is well defined (see Lemma \ref{Psi} above). It follows from
Lemmas \ref{P} and \ref{PHI} that for sufficiently large $\lambda>0$ and sufficiently small $r>0$, there holds $\Phi_\lambda(\sigma)\in
J_\rho^{-4L}$ for all $\sigma\in\Sigma_{\epsilon,k}$. This together with \eqref{betarho} implies
\be\label{beta}\beta_{\lambda,\rho}\leq -4L.\ee

Now we  claim $\alpha_{\lambda,\rho}>-2L$. For otherwise, $\alpha_{\lambda,\rho}\leq-2L$. By definition of $\alpha_{\lambda,\rho}$, namely
\eqref{min-max},
there exists some $\gamma_1\in\Gamma_\lambda$ such that
$\sup_{(\sigma,t)\in\widehat{\Sigma}_{\epsilon,k}}J_\rho(\gamma_1(\sigma,t))\leq-{3}L/2$. As a consequence
$$J_\rho(\gamma_1(\sigma,t))\leq -\f{3}{2}L\quad{\rm for\,\,all}\quad (\sigma,t)\in\widehat{\Sigma}_{\epsilon,k}.$$
Since $\Psi:J_\rho^{-L}\ra\overline\Sigma_k$ is continuous, the map $\Psi\circ\gamma_1:\widehat{\Sigma}_{\epsilon,k}\ra\overline\Sigma_{k}$ is also continuous.
Noting that $\gamma_1(\sigma,1)=\Phi_\lambda(\sigma)$ and $\gamma_1(\sigma,0)\equiv u_0\in J^{-L}_\rho$ for all $\sigma\in \Sigma_{\epsilon,k}$, if we let $\pi:\overline\Sigma_k\ra \Sigma_{\epsilon,k}$ be a continuous projection, then
$\pi\circ\Psi\circ\Phi_\lambda: \Sigma_{\epsilon,k}\ra \Sigma_{\epsilon,k}$ is homotopic to
a constant map $\pi\circ\Psi\circ\gamma_1(\cdot,0):\Sigma_{\epsilon,k}\ra\Sigma_{\epsilon,k}$. Moreover, by Lemma \ref{PHI},
$\pi\circ\Psi\circ\Phi_\lambda$ is homotopic to ${\rm Id}:\Sigma_{\epsilon,k}\ra \Sigma_{\epsilon,k}$. Hence
the identity map ${\rm Id}:\Sigma_{\epsilon,k}\ra \Sigma_{\epsilon,k}$ is homotopic to the constant map
$\pi\circ\Psi\circ\gamma_1(\cdot,0):\Sigma_{\epsilon,k}\ra \Sigma_{\epsilon,k}$, which contradicts the fact that $\Sigma_{\epsilon,k}$
is non-contractible.
Therefore
\be\label{alpha}\alpha_{\lambda,\rho}>-2L.\ee

Since $J_\rho\in C^2(W_0^{1,2}(\Sigma),\mathbb{R})$ and $\widehat{\Sigma}_{\epsilon,k}$ is a compact metric space, we immediately have that
$\beta_{\lambda,\rho}>-\infty$ and $\alpha_{\lambda,\rho}<+\infty$.
 This together with \eqref{beta} and \eqref{alpha} concludes the lemma. $\hfill\Box$\\

 To proceed, we need several uniform estimates for functionals $J_\rho$.

 \subsection{Uniform estimates with respect to \texorpdfstring{$\rho$}{rho}}\label{2.5}

 Let $[a,b]\subset(8k\pi,8(k+1)\pi)$ be any closed interval. Let $L>0$ be sufficiently large such that
 \be\label{psi-b}\Psi: J_b^{-L}\ra \overline\Sigma_k\ee
 is a continuous map defined as in Lemma \ref{Psi}. Let $\Sigma_{\epsilon,k}$ be given as in the previous section.
 Choose a sufficiently large $\lambda>0$  such that for all $\sigma\in \Sigma_{\epsilon,k}$, $\Phi_\lambda(\sigma)=\phi_{\lambda,\sigma}$ satisfies
  \be\label{Ja}J_a(\phi_{\lambda,\sigma})\leq -(a-8k\pi)\log\lambda\leq -4L,\ee
  where $\phi_{\lambda,\sigma}$ is defined as in \eqref{p-1}. It should be remarked that the choice of $\lambda$ depends not only on $L$, $k$, $a$, but also on
  $\epsilon$. Let $\Gamma_\lambda$ and $\alpha_{\lambda,\rho}$ be given as in \eqref{path} and \eqref{min-max} respectively.

 \begin{lemma}\label{uniform-1}
 Let $\rho\in[a,b]$. Then $\Psi:J_\rho^{-L}\ra \overline\Sigma_k$ is well defined uniformly with respect to $\rho$. Moreover, for all
 $\rho\in[a,b]$, there holds
 $$J_\rho(\phi_{\lambda,\sigma})\leq -4L,\quad\forall \sigma\in \Sigma_{\epsilon,k}.$$
 \end{lemma}
 \proof Let $\rho\in [a,b]$. If $u\in J_\rho^{-L}$, then $J_\rho(u)\leq -L$. This implies
 $$\log\int_\Sigma he^udv_g>0.$$
 It follows that $J_b(u)\leq J_\rho(u)\leq -L$, and that $u\in J_b^{-L}$. As a consequence $J_\rho^{-L}\subset J_b^{-L}$, and thus
 by \eqref{psi-b},
 $\Psi:J_\rho^{-L}\ra\overline\Sigma_k$ is well defined.

 Let $\sigma\in\Sigma_{\epsilon,k}$ and $\phi_{\lambda,\sigma}$ satisfy \eqref{Ja}. Similarly as above, we have by \eqref{Ja},
 \bna
 J_\rho(\phi_{\lambda,\sigma})&=&\f{1}{2}\int_\Sigma|\nabla_g\phi_{\lambda,\sigma}|^2dv_g-\rho\int_\Sigma he^{\phi_{\lambda,\sigma}}dv_g\\
 &\leq&\f{1}{2}\int_\Sigma|\nabla_g\phi_{\lambda,\sigma}|^2dv_g-a\int_\Sigma he^{\phi_{\lambda,\sigma}}dv_g\\
 &=&J_a(\phi_{\lambda,\sigma})\leq -4L.
 \ena
 This ends the proof of the lemma. $\hfill\Box$\\

 For simplicity, we denote $\alpha_{\lambda,\rho}$ by $\alpha_\rho$. By Lemmas \ref{mini} and \ref{uniform-1},
 $\alpha_\rho$ is a real number for any $\rho\in [a,b]$. Then
 we have an analog of \cite[Lemma 2.4]{DJLW2}, namely

 \begin{lemma}\label{monotonicity}
 $\alpha_\rho/\rho$ is decreasing in $\rho\in [a,b]$.
 \end{lemma}
 \proof Let $a\leq \rho_1<\rho_2\leq b$. Then for any $(\sigma,t)\in\widehat{\Sigma}_{\epsilon,k}$ and any $\gamma\in\Gamma_\lambda$,
 there holds
 $$
 \f{J_{\rho_1}(\gamma(\sigma,t))}{\rho_1}-\f{J_{\rho_2}(\gamma(\sigma,t))}{\rho_2}=\le(\f{1}{\rho_1}-\f{1}{\rho_2}\ri)\int_\Sigma
 |\nabla_g\gamma(\sigma,t)|dv_g\geq 0.
 $$
 It then follows that $\alpha_{\rho_1}/\rho_1\geq \alpha_{\rho_2}/\rho_2$. $\hfill\Box$\\

 By Lemma \ref{monotonicity}, $\alpha_\rho/\rho$ is differentiable almost everywhere in $[a,b]\subset(8k\pi,8(k+1)\pi)$. Denote
 \be\label{lambda}\Lambda_{a,b}=\le\{\rho\in(a,b): \f{\alpha_\rho}{\rho}\,\,{\rm is\,\, differentiable\,\, at}\,\,\rho\ri\}.\ee
 Then $\Lambda_{a,b}$ is a dense subset of $[a,b]$.

 \subsection{Existence for a dense set}\label{2.6}

 In this subsection, we shall prove that $J_\rho$ has a critical point for any $\rho\in \Lambda_{a,b}$. The argument we shall use
 is adapted from  Ding-Jost-Li-Wang \cite{DJLW2}. For readers' convenience, we provide the details here.
 \begin{lemma}\label{bdd}
 If $\rho\in\Lambda_{a,b}$, then $\alpha_\rho$ is differentiable at $\rho$. In particular, if $\rho\in\Lambda_{a,b}$, then we have
 $$\alpha_{\tilde{\rho}}=\alpha_\rho+O(\tilde{\rho}-\rho) \quad {\rm as}\,\,\, \tilde{\rho}\ra \rho.$$
 \end{lemma}
 \proof In view of \eqref{lambda}, it suffices to notice that $\alpha_\rho=\rho(\alpha_\rho/\rho)$. $\hfill\Box$\\

 As an analog of \cite[Lemma 3.2]{DJLW2}, we have the following:

 \begin{lemma}\label{dense}
 If $\rho\in\Lambda_{a,b}$, then $\alpha_\rho$ is a critical value of $J_\rho$.
 \end{lemma}
 \proof Let $(\rho_n)\subset [a,b]$ be an increasing sequence converging to $\rho\in\Lambda_{a,b}$. By the definition of $\alpha_{\rho_n}$,
 there must be a path $\gamma_n\in\Gamma_\lambda$ such that
 \be\label{rho-n}\sup_{u\in\gamma_n(\widehat{\Sigma}_{\epsilon,k})}J_{\rho_n}(u)\leq \alpha_{\rho_n}+\rho-\rho_n.\ee
 Also we use the definition of $\alpha_\rho$ to find some $u_n\in \gamma_n(\widehat{\Sigma}_{\epsilon,k})\subset W_0^{1,2}(\Sigma)$ with
 \be\label{J-geq}J_\rho(u_n)\geq \alpha_\rho-(\rho-\rho_n).\ee
  For the above $u_n$, we have by Lemma \ref{bdd},
  \bea\nonumber
  \f{1}{2}\int_\Sigma|\nabla_gu_n|^2dv_g&=&\f{\f{J_{\rho_n}(u_n)}{\rho_n}-\f{J_{\rho}(u_n)}{\rho}}{\f{1}{\rho_n}-\f{1}{\rho}}\\
  &\leq& \rho\f{\alpha_{\rho_n}-\alpha_\rho}{\rho-\rho_n}+\le(\f{1}{\rho_n}-\f{1}{\rho}\ri){\alpha_\rho}+\rho_n+\rho\nonumber\\
  &\leq&c_0\label{bbb}
  \eea
  for some constant $c_0$ depending only on $\rho$, $\alpha_\rho$ and $(\alpha_\rho/\rho)^\prime$. Moreover,
  by Lemmas \ref{monotonicity} and \ref{bdd}, and the estimate \eqref{rho-n}, one finds
  \be\label{ll}
  J_\rho(u_n)\leq\f{\rho}{\rho_n}J_{\rho_n}(u_n)
  \leq\f{\rho}{\rho_n}(\alpha_{\rho_n}+\rho-\rho_n)
  \leq\alpha_\rho+C(\rho-\rho_n)
  \ee
  for some constant $C$ independent of $n$.

  Suppose $\alpha_\rho$ is not a critical value of $J_\rho$. Since any bounded Palais-Smale sequence  must converge
  to a critical point of $J_\rho$ (see \cite[Lemma 3.1]{DJLW2}), there would exist $\delta>0$ such that
  \be\label{non-0}\|dJ_\rho(u)\|_{(W_0^{1,2}(\Sigma))^\ast}\geq 2\delta\ee
  for all $u\in\mathscr{N}_\delta$, where
  \be\label{Ndelta}\mathscr{N}_\delta=\le\{u\in W_0^{1,2}(\Sigma):\int_\Sigma|\nabla_gu|^2dv_g\leq 2c_0,\,|J_\rho(u)-\alpha_\rho|<\delta\ri\}.\ee
  It follows from \eqref{J-geq}, \eqref{bbb} and \eqref{ll} that $\mathscr{N}_\delta\not=\varnothing$. Let $\mathbf{X}_\rho:
  \mathscr{N}_\delta\ra W_0^{1,2}(\Sigma)$ be a pseudo-gradient vector field for $J_\rho$ in $\mathscr{N}_\delta$, namely a locally Lipschitz vector field satisfying $\|\mathbf{X}_\rho\|_{W_0^{1,2}(\Sigma)}\leq 1$ and
  \be\label{-delta}dJ_\rho(u)\le(\mathbf{X}_\rho(u)\ri)\leq -\delta.\ee
  Here we have used \eqref{non-0}. One can check that as $n\ra\infty$, $dJ_{\rho_n}(u)$ converges to $dJ_\rho(u)$ in $(W_0^{1,2}(\Sigma))^\ast$
  uniformly in $u$ with $\int_\Sigma |\nabla_gu|^2dv_g\leq c^\ast$. Thus $\mathbf{X}_\rho$ is also a pseudo-gradient vector
  field for $J_{\rho_n}$ in $\mathscr{N}_\delta$. Moreover, there holds for all $u\in\mathscr{N}_\delta$ and sufficiently large $n$,
  \be\label{dj-n}dJ_{\rho_n}(u)\le(\mathbf{X}_\rho(u)\ri)\leq -\delta/2.\ee
  Take a Lipschitz continuous cut-off function $\eta$ such that $0\leq\eta\leq 1$, $\eta\equiv 1$ in $\mathscr{N}_{\delta/2}$,
  and $\eta\equiv 0$ outside $\mathscr{N}_{\delta}$. Let $\psi:W_0^{1,2}(\Sigma)\times [0,+\infty)$ be the flow generated by
  $\eta\mathbf{X}_\rho$, which reads as
  $$\le\{\begin{array}{lll}
  \f{\p}{\p s}\psi(u,s)=\eta(\psi(u,s))\mathbf{X}_\rho(\psi(u,s))\\[1.5ex]
  \psi(u,0)=u.
  \end{array}\ri.$$
  This flow has long time existence because it always remains stationary outside $\mathscr{N}_\delta$. It follows from
  \eqref{-delta} that for all $u\in\mathscr{N}_{\delta/2}$,
  \be\label{deriv-0}\le.\f{d}{ds}\ri|_{s=0}J_\rho(\psi(u,s))=dJ_\rho(u)\le(\mathbf{X}_\rho(u)\ri)\leq -\delta.\ee

  In view of \eqref{alpha} and \eqref{Ndelta}, one easily sees $\phi_{\lambda,\sigma}\not\in \mathscr{N}_\delta$ for all
  $\sigma\in \Sigma_{\epsilon,k}$. Since $J_\rho(\psi(\phi_{\lambda,\sigma},s))$ is decreasing in $s$, there holds
  $J_\rho(\psi(\phi_{\lambda,\sigma},s))\leq -4L$ for all $s\in[0,+\infty)$. Hence $\psi(\gamma_n(\sigma,1),s)\not\in\mathscr{N}_\delta$,
  and thus $\psi(\gamma_n(\sigma,1),s)\equiv \psi(\gamma_n(\sigma,1),0)=\phi_{\lambda,\sigma}$ for all $\sigma\in \Sigma_{\epsilon,k}$ and all $s\in[0,+\infty)$. As a consequence, if we write $\psi_s(\cdot)=\psi(\cdot,s)$, then we have $\psi_s\circ\gamma_n\in\Gamma_\lambda$.
  By \eqref{ll} and the monotonicity of $J_\rho(\psi_s(u))$ in $s$, we obtain
  \be\label{bound}
  \alpha_\rho\leq\sup_{u\in\psi_s\circ\gamma_n(\widehat{\Sigma}_{\epsilon,k})}J_\rho(u)\leq
  \sup_{u\in\gamma_n(\widehat{\Sigma}_{\epsilon,k})}J_\rho(u)\leq \alpha_\rho+C(\rho-\rho_n).
  \ee

  We now claim that
 \be\label{usn}\sup_{u\in\psi_s\circ\gamma_n(\widehat{\Sigma}_{\epsilon,k})}J_\rho(u)\,\, {\rm is\,\, achieved\,\, in}\,\, \mathscr{N}_{\delta/2}.\ee
  In fact, since $\widehat{\Sigma}_{\epsilon,k}$ is a compact metric space, the continuous function
  $J_\rho(\psi_s\circ\gamma_n(\cdot,\cdot))$ attains its supremum at some $(\sigma_{0},t_{0})\in\widehat{\Sigma}_{\epsilon,k}$.
  As a result, the function $u_{s,n}=\psi_s\circ\gamma_n(\sigma_0,t_0)$ achieves
  $\sup_{u\in\psi_s\circ\gamma_n(\widehat{\Sigma}_{\epsilon,k})}J_\rho(u)$.
    If $n$ is chosen sufficiently large, \eqref{bound} implies that
  $\alpha_\rho\leq J_\rho(u_{s,n})\leq \alpha_\rho+\delta/2$. By \eqref{dj-n}, $J_{\rho_n}(\psi_s(u))$ is decreasing in $s$,
  which together with \eqref{rho-n} gives $J_{\rho_n}(u_{s,n})\leq \alpha_{\rho_n}+\rho-\rho_n$. It then follows that
  \bna
  \f{1}{2}\int_\Sigma |\nabla_gu_{s,n}|^2dv_g&=&\le(\f{1}{\rho_n}-\f{1}{\rho}\ri)^{-1}\le(\f{J_{\rho_n}(u_{s,n})}{\rho_n}-
  \f{J_{\rho}(u_{s,n})}{\rho}\ri)\\
  &\leq&\le(\f{1}{\rho_n}-\f{1}{\rho}\ri)^{-1}\le(\f{\alpha_{\rho_n}+\rho-\rho_n}{\rho_n}-\f{\alpha_\rho}{\rho}\ri)\\
  &\leq&c_0,
  \ena
  where $c_0$ is the same constant as in \eqref{bbb}. Therefore $u_{s,n}\in\mathscr{N}_{\delta/2}$, and
  our claim is confirmed.

  Let $\bar{s}>0$ and $\overline{\gamma}_n=\psi_{\bar{s}}\circ\gamma_n$. Then  we have by \eqref{deriv-0} and
  \eqref{usn} that
  \bea
  \le.\f{d}{ds}\ri|_{s=\bar{s}}\sup_{(\sigma,t)\in\widehat{\Sigma}_{\epsilon,k}}J_\rho(\psi_s\circ\gamma_n(\sigma,t))
  &=&\le.\f{d}{ds}\ri|_{s=\bar{s}}\sup_{(\sigma,t)\in\widehat{\Sigma}_{\epsilon,k}}J_\rho(\psi_{s-\bar{s}}\circ\overline\gamma_n(\sigma,t))\nonumber\\
  &=&\le.\f{d}{ds}\ri|_{s=0}\sup_{(\sigma,t)\in\widehat{\Sigma}_{\epsilon,k}}J_\rho(\psi_{s}\circ\overline\gamma_n(\sigma,t))\nonumber\\
  &\leq&\sup_{u\in\mathscr{N}_{\delta/2}}\le.\f{d}{ds}\ri|_{s=0}J_\rho(\psi_s(u))\nonumber\\
  &\leq&-\delta.\label{fin}
  \eea
  Using the Newton-Lebnitz formula, we conclude from \eqref{bound} and \eqref{fin} that
  $$\sup_{u\in\psi_s\circ\gamma_n(\widehat{\Sigma}_{\epsilon,k})}J_\rho(u)<\alpha_\rho,$$
  if  $s>0$ is sufficiently large. This contradicts the definition of $\alpha_\rho$,
  and ends the proof of the lemma. $\hfill\Box$

   \subsection{Existence for all \texorpdfstring{$\rho\in(8k\pi,8(k+1)\pi)$}{non critical cases}}\label{2.7}

   In this subsection, we use previous analysis to complete the proof of \autoref{thm1}. \\

   \noindent{\it Proof of \autoref{thm1}.}
  For any $\rho\in (8k\pi,8(k+1)\pi)$, $k\in\mathbb{N}^\ast$, there are two constants $a$ and $b$ with $8k\pi<a<\rho<b<8(k+1)\pi$. By Lemma \ref{dense}, we may take an increasing sequence of numbers
  $(\rho_n)\subset\Lambda_{a,b}$ such that  $\rho_n\ra \rho$ and $\alpha_{\rho_n}$ is achieved by $u_n\in W_0^{1,2}(\Sigma)$. Moreover,
  $u_n$ satisfies
  \be\label{n-equ}\Delta_gu_n=\rho_n\f{he^{u_n}}{\int_\Sigma he^{u_n}dv_g}\quad{\rm in}\quad\Sigma.\ee
  By Lemma \ref{monotonicity},
  \be\label{alpha-n}\alpha_{\rho_n}\leq \f{b}{a}\alpha_a.\ee
  Denoting $v_n=u_n-\log\int_\Sigma he^{u_n}dv_g$, we have by \eqref{n-equ} that
  $$\le\{\begin{array}{lll}
  \Delta_gv_n=\rho_nhe^{v_n}\\[1.5ex]
  \int_\Sigma he^{v_n}dv_g=1.
  \end{array}\ri.$$
  By Lemma \ref{compact}, $(u_n)$ is bounded in $L^\infty(\overline{\Sigma})$. Let $\Omega_1,\cdots,\Omega_{k+1}$ be disjoint sub-domains
  of $\overline\Sigma$. By Lemma \ref{improved},
  $$\log\int_\Sigma e^{u_n}dv_g\leq \f{1}{16(k+1)\pi-\epsilon}\int_\Sigma|\nabla_gu_n|^2dv_g+C_\epsilon$$
  for any $\epsilon>0$ and some constant $C_\epsilon>0$. This together with (\ref{alpha-n}) implies that for $0<\epsilon<16(k+1)\pi-2b$,
  \bna
  \f{1}{2}\int_\Sigma|\nabla_gu_n|^2dv_g&=&J_{\rho_n}(u_n)+\rho_n\log\int_\Sigma he^{u_n}dv_g\\
  &\leq&\f{b}{16(k+1)\pi-\epsilon}\int_\Sigma|\nabla_gu_n|^2dv_g+C.
  \ena
  It then follows that $(u_n)$ is bounded in $W_0^{1,2}(\Sigma)$. With no loss of generality, we assume
  $u_n$ converges to $u_\rho$ weakly in $W_0^{1,2}(\Sigma)$, strongly in $L^p(\Sigma)$ for any $p>1$, and almost everywhere in $\Sigma$.
  Moreover, $e^{u_n}$ converges to $e^{u_\rho}$ strongly in $L^p(\Sigma)$ for  any $p>1$. By \eqref{n-equ}, $u_\rho$ is a distributional
  solution of \eqref{Dir}.
  Hence $u_\rho$ is a critical point of $J_\rho$. $\hfill\Box$

\section{The Neumann boundary value problem}\label{Neuma}
In this section, we shall prove \autoref{thm2} by the min-max method. Since part of the proof is analogous to that of \autoref{thm1}, we only give its outline but stress the difference. In \autoref{3.1}, we prove a compactness result for solutions of \eqref{num}; In  \autoref{3.2}, we derive an improved Trudinger-Moser inequality for functions $u\in W^{1,2}(\Sigma)$ with $\int_\Sigma udv_g=0$; In  \autoref{3.3}, we construct two continuous maps between sub-levels of $J_\rho$ and the topological space $\mathscr{S}_k$, where $J_\rho$ and $\mathscr{S}_k$ are defined as in \eqref{funct-2} and \eqref{subset} respectively; In \autoref{3.4}, we construct min-max levels of $J_\rho$; The remaining part of the proof of \autoref{thm2} is outlined in \autoref{3.5}.

\subsection{Compactness analysis}\label{3.1}
Let $(\rho_n)$ be a number sequence tending to $\rho\in\mathbb{R}$, $(h_n)$ be a function sequence converging to
$h$ in $C^1(\overline{\Sigma})$, and $(u_n)$ be a sequence of solutions to
\be\label{Neumann}\le\{\begin{array}{lll}
\Delta_g u_n=\rho_n\le(\f{h_ne^{u_n}}{\int_\Sigma h_ne^{u_n}dv_g}-\f{1}{|\Sigma|}\ri)&{\rm in}&\Sigma\\[1.5ex]
\p u_n/\p\mathbf{v}=0&{\rm on}&\p\Sigma\\[1.5ex]
\int_\Sigma u_ndv_g=0.
\end{array}\ri.\ee
Denote $v_n=u_n-\log\int_\Sigma h_ne^{u_n}dv_g$. Then
$\Delta_{g} v_n=\rho_n(h_ne^{v_n}-1/|\Sigma|)$ and $\int_\Sigma h_ne^{v_n}dv_g=1$. Concerning the compactness of $(u_n)$, we have an analog of
Lemma \ref{compact}, namely

\begin{lemma}\label{compact-2}
Assume $\rho$ is a positive number and $h$ is a positive function. Up to a subsequence, there holds one of the following alternatives:\\
$(i)$ $(u_n)$ is bounded in $L^\infty(\overline\Sigma)$;\\
$(ii)$ $(v_n)$ converges to $-\infty$ uniformly in $\overline\Sigma$;\\
$(iii)$ there exists a finite singular set $\mathcal{S}=\{p_1,\cdots,p_m\}\subset {\overline\Sigma}$ such that
for any $1\leq j\leq m$, there is a sequence of points $\{p_{j,n}\}\subset\overline\Sigma$ satisfying
$p_{j,n}\ra p_j$, $u_n(p_{j,n})\ra+\infty$, and $v_n$ converges to $-\infty$ uniformly on any compact subset of
$\overline{\Sigma}\setminus\mathcal{S}$ as $n\ra\infty$. Moreover, if $\mathcal{S}$ has $\ell$ points in $\Sigma$
and $(m-\ell)$ points on $\p\Sigma$, then
$$\rho_n\int_\Sigma h_ne^{v_n}dv_g\ra 4(m+\ell)\pi. $$
\end{lemma}

\proof We modify arguments in the proof of Lemma \ref{compact}, and divide it into two parts.\\

{\bf Part I. Analysis in the interior domain $\Sigma$}

Let $(u_n)$ be a sequence of solutions to \eqref{Neumann}. By the Green representation formula (see \cite{Yang-Zhou}), we have
\be\label{W1q-2}\|u_n\|_{W^{1,q}(\Sigma)}\leq C_q,\quad\forall 1<q<2.\ee
Since $(h_n)$ converges to $h>0$ in $C^1(\overline\Sigma)$, there exists some constant $C>0$ such that
for all $n\in \mathbb{N}$,
\be\label{lower-b-2}\int_\Sigma e^{v_n}dv_g\leq C.\ee
Moreover, Jensen's inequality implies 
$$\liminf_{n\ra\infty}\int_\Sigma e^{u_n}dv_g\geq|\Sigma|.$$

We assume with no loss of generality, $\rho_nh_ne^{v_n}dv_g$ converges to some nonnegative measure $\mu$ on $\overline\Sigma$. If $\mu(x^\ast)<4\pi$ for some $x^\ast\in\Sigma$,
then there exist two positive constants $\epsilon_0$ and $r_0$ verifying
$$\int_{B_{x^\ast}(r_0)}\rho_nh_ne^{v_n}dv_g\leq 4\pi-\epsilon_0.$$
In view of \eqref{Neumann}, by a result of Brezis-Merle \cite[Theorem 1]{B-M} and elliptic estimates, we have that
$(u_n)$ is bounded in $L^{\infty}(B_{x^\ast}(r_0/2))$. This leads to $\mu(x^\ast)=0$. Define a set
$\mathcal{S}=\{x\in\Sigma: \mu(x)\geq 4\pi\}$.
If $\mathcal{S}\not=\varnothing$, then by almost the same argument as the proof of \eqref{local-uniform}, we conclude
that for any compact set $A\subset\Sigma\setminus\mathcal{S}$, there
holds
\be\label{local-uniform-2}v_n\ra-\infty\,\,{\rm uniforlmy\,\,in}\,\, x\in A.\ee

  Assume  $\mathcal{S}=\{x_1,\cdots,x_j\}$ for some positive integer $j$. We shall show that $\mu(x_i)=8\pi$ for all $1\leq i\leq j$.
With no loss of generality, it suffices to prove $\mu(x_1)=8\pi$. For this purpose, we choose an isothermal coordinate system
$\phi:U\ra\mathbb{B}_1=\{(y_1,y_2)\in\mathbb{R}^2: y_1^2+y_2^2<1\}$ near $x_1$. In such coordinates, the metric $g$ and the Laplace-Beltrami operator $\Delta_g$
are represented by $g=e^{\psi(y)}(dy_1^2+dy_2^2)$ and
 $\Delta_g=-e^{-\psi(y)}\Delta_{\mathbb{R}^2}$ respectively, where $\psi$ is a smooth function with $\psi(0,0)=0$, and $\Delta_{\mathbb{R}^2}={\p^2}/{\p y_1^2}+{\p^2}/{\p y_2^2}$ denotes the standard Laplacian on $\mathbb{R}^2$.
 Set $\widetilde{u}=u\circ \phi^{-1}$ for any function $u:U\ra\mathbb{R}$. Since $(u_n)$ satisfies \eqref{Neumann}, then
 $\widetilde{u}_n=u_n\circ \phi^{-1}$ satisfies
 \be\label{local-2}-\Delta_{\mathbb{R}^2}\widetilde{u}_n(y)=e^{\psi(y)}\rho_n(\widetilde{h}_n(y)e^{\widetilde{v}_n(y)}-|\Sigma|^{-1}),
 \quad y\in\mathbb{B}_1.\ee
 Multiplying both sides of \eqref{local-2} by $y\cdot\nabla_{\mathbb{R}^2}\widetilde{u}_n(y)$, we have by integration by parts 
 \bea\label{Pohozaev-2}\nonumber\f{r}{2}\int_{\p\mathbb{B}_r}|\nabla_{\mathbb{R}^2}\widetilde{u}_n|^2d\sigma-r\int_{\p\mathbb{B}_r}
 \langle\nabla_{\mathbb{R}^2}\widetilde{u}_n,\mathbf{\nu} \rangle^2d\sigma&=&r\int_{\p\mathbb{B}_r}
 e^{\psi}\rho_n\widetilde{h}_ne^{\widetilde{v}_n}d\sigma-\int_{\mathbb{B}_r}e^{\widetilde{v}_n}
 \rho_n\langle\nabla_{\mathbb{R}^2}(e^{\psi}\widetilde{h}_n),y \rangle dy\\&&-2\int_{\mathbb{B}_r}
 e^{\psi}\rho_n\widetilde{h}_ne^{\widetilde{v}_n}dy+\f{\rho_n}{|\Sigma|}\int_{\mathbb{B}_r}e^{\psi}y\cdot\nabla_{\mathbb{R}^2}\widetilde{u}_ndy,\eea
 where $\mathbb{B}_r=\{y\in\mathbb{R}^2:y_1^2+y_2^2<r\}$, $\p\mathbb{B}_r=\{y\in\mathbb{R}^2:y_1^2+y_2^2=r\}$, and
 $\nu$ denotes the unit outward vector on $\p\mathbb{B}_r$. In view of \eqref{local-uniform-2}, $(u_n)$ converges to
 a Green function $G(x_1,\cdot)$ weakly in $W^{1,q}(\Sigma)$, and in $C^2_{\rm loc}(\Sigma\setminus\mathcal{S})$. Locally
 $G(x_1,\cdot)$ satisfies
 $$\Delta_{g,z}G(x_1,z)=\mu(x_1)\delta_{x_1}(z)-\rho|\Sigma|^{-1},\quad\forall z\in \phi^{-1}(\mathbb{B}_{1}).$$
 Clearly $\widetilde{G}(y)=G(x_1,\phi^{-1}(y))=-\f{\mu(x_1)}{2\pi}\log|y|+\eta(y)$ for some $\eta\in C^2(\mathbb{B}_1)$.
 Passing to the limit $n\ra\infty$ first, and then $r\ra 0$ in \eqref{Pohozaev-2}, we obtain
 \be\label{mess-2}\mu(x_1)=\lim_{r\ra 0}\le(\f{r}{2}\int_{\p\mathbb{B}_r}\langle\nabla_{\mathbb{R}^2}\widetilde{G},\nu\rangle^2d\sigma-\f{r}{4}
 \int_{\p\mathbb{B}_r}|\nabla_{\mathbb{R}^2}\widetilde{G}|^2d\sigma\ri)=\f{(\mu(x_1))^2}{8\pi}.\ee
 This immediately leads to $\mu(x_1)=8\pi$. In conclusion, we have
 \be\label{inner-2}\mu(x_i)=8\pi\,\,\,{\rm for\,\,all}\,\,\, 1\leq i\leq j.\ee

{\bf Part II. Analysis on the boundary $\p\Sigma$}

 Let $x^\ast\in \overline\Sigma$ be fixed. Note that $\rho_nh_ne^{v_n}dv_{g}$ converges to the nonnegative Radon measure $\mu$ on $\overline{\Sigma}$
 as $n\ra\infty$.
 If $\mu(x^\ast)<2\pi$, there exist a neighborhood $V$ of $x^\ast$ and a number $\gamma_0>0$ such that
 \be\label{2pi-2}\int_V\rho_nh_ne^{v_n}dv_g\leq 2\pi-\gamma_0.\ee
 With no loss of generality, we take an isothermal coordinate system $(V,\phi,\{y_1,y_2\})$ such that
 $\phi(x^\ast)=(0,0)$, and $\phi:V\ra{\mathbb{B}_1^+}\cup\Gamma=\{(y_1,y_2):y_1^2+y_2^2< 1,\,y_2\geq 0\}$, where
 $\Gamma=\{(y_1,y_2): |y_1|<1, y_2=0\}$. Moreover, in this coordinate system, the metric $g=e^{\psi(y)}(dy_1^2+dy_2^2)$, and the Laplace-Beltrami operator $\Delta_g=-e^{-\psi(y)}\Delta_{\mathbb{R}^2}$, where $\psi:\mathbb{B}_1^+\cup\Gamma\ra \mathbb{R}$ is a smooth function
 with $\psi(0,0)=0$; moreover, $\p/\p\mathbf{v}=e^{-\psi(y)/2}\p/\p y_2$. For more details about isothermal coordinates on the boundary, we refer the readers to \cite[Section 2]{Yang-Zhou}. Now the local version of \eqref{Neumann} reads as
 \be\label{semi-equation-2}\le\{\begin{array}{lll}
 -\Delta_{\mathbb{R}^2}(u_n\circ\phi^{-1})(y)=e^{\psi(y)}\rho_n\le((h_n\circ\phi^{-1})(y)e^{(v_n\circ\phi^{-1})(y)}-|\Sigma|^{-1}\ri)
 &{\rm in}&\mathbb{B}_1^+\\[1.5ex]
 \f{\p}{\p y_2}(u_n\circ\phi^{-1})(y)=0&{\rm on}& \Gamma.
 \end{array}\ri.\ee
 For any function $u:V\ra \mathbb{R}$, we define a function $\widetilde{u}: \mathbb{B}_1\ra\mathbb{R}$ by
 \be\label{u-tilde-2}\widetilde{u}(y_1,y_2)=\le\{
 \begin{array}{lll}
 u\circ\phi^{-1}(y_1,y_2)&{\rm if}& y_2\geq 0\\[1.5ex]
 u\circ\phi^{-1}(y_1,-y_2)&{\rm if}& y_2< 0.
 \end{array}\ri.\ee
 One can easily derive from \eqref{semi-equation-2} that $\widetilde{u}_n$ is a distributional solution of
 \be\label{un-2}-\Delta_{\mathbb{R}^2}\widetilde{u}_n(y)=\widetilde{f}_n(y),\quad y\in\mathbb{B}_1,\ee
 where $\widetilde{f}_n$ is defined as in \eqref{u-tilde-2} and for $y\in \mathbb{B}_1^+\cup\Gamma$,
 $$f_n\circ\phi^{-1}(y)=e^{\psi(y)}\rho_n\le((h_n\circ\phi^{-1})(y)e^{(v_n\circ\phi^{-1})(y)}-|\Sigma|^{-1}\ri).$$
 In view of \eqref{2pi-2} and the fact $\psi(0,0)=0$, there would exist a number $0<r_0<1$ such that
 $$\int_{\mathbb{B}_{r_0}}|\widetilde{f}_n(y)|dy\leq 4\pi-\gamma_0.$$
 Let $w_n$ be a solution of
 $$\le\{\begin{array}{lll}
 -\Delta_{\mathbb{R}^2}w_n=\widetilde{f}_n&{\rm in}&\mathbb{B}_{r_0}\\[1.5ex]
 w_n=0&{\rm on}&\p\mathbb{B}_{r_0}.
 \end{array}\ri.$$
 By \cite[Theorem 1]{B-M}, there exists some constant $C$ depending only on $\epsilon_0$
 and $r_0$ such that
 $$\int_{\mathbb{B}_{r_0}}\exp\le(\f{(4\pi-\gamma_0/2)|w_n|}{\|\widetilde{f}_n\|_{L^1(\mathbb{B}_{r_0})}}\ri)dy\leq C.$$
 Hence there exists some $q_0>1$ such that
 \be\label{Lq0-2}\|e^{|w_n|}\|_{L^{q_0}(\mathbb{B}_{r_0})}\leq C.\ee
 Let
 $\eta_n=\widetilde{u}_n-w_n$. Then $\eta_n$ satisfies
 \be\label{harmonic-2}\le\{\begin{array}{lll}
 -\Delta_{\mathbb{R}^2}\eta_n=0&{\rm in}&\mathbb{B}_{r_0}\\[1.5ex]
 \eta_n=\widetilde{u}_n&{\rm on}&\p\mathbb{B}_{r_0}.
 \end{array}\ri.\ee
 Noticing \eqref{W1q-2} and \eqref{Lq0-2}, we have by applying elliptic estimates to \eqref{harmonic-2} that
 \be\label{etan-2}\|\eta_n\|_{L^\infty(\mathbb{B}_{r_0/2})}\leq C.\ee
 Combining \eqref{lower-b-2}, \eqref{Lq0-2} and \eqref{etan-2}, we conclude
 $\|\widetilde{f}_n\|_{L^{q_0}(\mathbb{B}_{r_0/2})}\leq C$. Applying elliptic estimates to \eqref{un-2}, we obtain
 $\|\widetilde{u}_n\|_{L^\infty(\mathbb{B}_{r_0/4})}\leq C$, which implies
 $\|{u}_n\|_{L^\infty(\phi^{-1}(\mathbb{B}^{+}_{r_0/4}))}\leq C$. In conclusion, we have that if $\mu(x^\ast)<2\pi$, then
 $(u_n)$ is uniformly bounded near $x^\ast$. This also leads to $\mu(x^\ast)=0$.

 If  $\mu(x^\ast)\geq 2\pi$, in the same coordinate system $(V,\phi,\{y_1,y_2\})$ as above, $\widetilde{f}_n(y)dy$ converges to a Radon measure $\widetilde{\mu}$ with
 $\widetilde{\mu}(0,0)=2\mu(x^\ast)\geq 4\pi$. Obviously there exists some $r_1>0$ such that for any $x\in \mathbb{B}_{r_1}\setminus\{(0,0)\}$,
 $\widetilde{\mu}(x)=0$. Using the same argument as the proof  of \eqref{local-uniform-2}, we conclude that
 for any compact set $A\subset \mathbb{B}_{r_1}\setminus\{(0,0)\}$, $\widetilde{v}_n$ converges to $-\infty$ uniformly in $A$. This leads to
 $\widetilde{f}_n(y)dy$ converges to the Dirac measure $\widetilde{\mu}(0,0)\delta_{(0,0)}(y)$. Recalling \eqref{W1q-2}, we have
 $\widetilde{u}_n$ converges to $\widetilde{G}_0$ weakly in $W^{1,q}(\mathbb{B}_{r_1})$ and a.e. in $\mathbb{B}_{r_1}$, where $\widetilde{G}_0$
 satisfies
 $$-\Delta_{\mathbb{R}^2}\widetilde{G}_0(y)=\widetilde{\mu}(0,0)\delta_{(0,0)}(y)-\rho|\Sigma|^{-1},\quad y\in\mathbb{B}_{r_1}.$$
 Clearly $\widetilde{G}_0$ is represented by
 $$\widetilde{G}_0(y)=-\f{\widetilde{\mu}(0,0)}{2\pi}\log|y|+A_0+O(|y|)$$
 as $y\ra 0$, where $A_0$ is a constant. Noting that $\widetilde{v}_n$ converges to $-\infty$
 locally uniformly in $\mathbb{B}_{r_1}\setminus\{(0,0)\}$, we have by applying elliptic estimates to \eqref{un-2} that
 $$\widetilde{u}_n\ra \widetilde{G}_0\quad {\rm in}\quad C^2_{\rm loc}(\mathbb{B}_{r_1}\setminus\{(0,0)\}).$$
 Multiplying both sides of \eqref{un-2} by $y\cdot\nabla_{\mathbb{R}^2}\widetilde{u}_n(y)$, completely analogous to \eqref{Pohozaev-2}
 and \eqref{mess-2}, we obtain
 $\widetilde{\mu}(0,0)=8\pi$, and thus
 \be\label{boundary-2}\mu(x^\ast)=4\pi.\ee

 Note that if $\mu(x_i)>0$ for some $x_i\in\overline\Sigma$, then there must exist $x_{i,n}\subset\overline\Sigma$ satisfying
 $u_n(x_{i,n})\ra+\infty$. For otherwise, $(u_n)$ is uniformly bounded near $x$, which leads to $\mu(x)=0$. The lemma then follows
 from \eqref{inner-2} and \eqref{boundary-2} immediately. $\hfill\Box$

 \subsection{An improved Trudinger-Moser inequality}\label{3.2}
 For a compact surface with smooth boundary, it was proved by Yang \cite{2006} that
 \be\label{T-M-mean-2}\sup_{u\in W^{1,2}(\Sigma),\,\int_\Sigma|\nabla_gu|^2dv_g\leq 1,\,\int_\Sigma udv_g=0}\int_\Sigma e^{2\pi u^2}dv_g<\infty.\ee
 Denote $\overline{u}=\f{1}{|\Sigma|}\int_\Sigma udv_g$. By \eqref{T-M-mean-2} and the Young inequality, we obtain
 \bna\log\int_\Sigma e^{u-\overline{u}}dv_g&\leq&\log\int_\Sigma e^{2\pi\f{(u-\overline{u})^2}{\|\nabla_gu\|_2^2}+\f{1}{8\pi}
 \|\nabla_gu\|_2^2}dv_g\\
 &=&\f{1}{8\pi}\int_\Sigma|\nabla_gu|^2dv_g+\log\int_\Sigma e^{2\pi\f{(u-\overline{u})^2}{\|\nabla_gu\|_2^2}}dv_g\\
 &\leq& \f{1}{8\pi}\int_\Sigma|\nabla_gu|^2dv_g+C.\ena
Hence
\be\label{T-M-weak-2}\log\int_\Sigma e^{u}dv_g\leq \f{1}{8\pi}\int_\Sigma|\nabla_gu|^2dv_g+\f{1}{|\Sigma|}\int_\Sigma udv_g,\quad\forall
u\in W^{1,2}(\Sigma).\ee

 \begin{lemma}\label{improved-2}
 Let $b_0$ and $\gamma_0$ be two positive constants, $\Omega_1,\cdots,\Omega_k$ be $k$ domains of $\overline\Sigma$ with
 ${\rm dist}(\Omega_i,\Omega_j)\geq b_0$ for all $1\leq i<j\leq k$. Then for any $\epsilon>0$, there exists some constant $C$
  depending only on $b_0,\gamma_0,k,\epsilon$, such that
 \be\label{imp-2}\log\int_\Sigma e^udv_g\leq \f{1}{8k\pi-\epsilon}\int_\Sigma|\nabla_gu|^2dv_g+\f{1}{|\Sigma|}\int_\Sigma udv_g+C\ee
 for all $u\in W^{1,2}(\Sigma)$ with
 \be\label{comp-2}\int_{\Omega_i}e^{u}dv_g\geq\gamma_0\int_\Sigma e^udv_g,\,\, i=1,\cdots,k.\ee
 \end{lemma}

 \proof  We follow the lines of Chen-Li \cite{Chen-Li}. Take smooth functions $\phi_1,\cdots,\phi_k$
 defined on $\overline\Sigma$ satisfying
 \be\label{phi-1-2}{\rm supp}\phi_i\cap {\rm supp}\phi_j=\varnothing,\quad\forall 1\leq i<j\leq k,\ee
 \be\label{phi-2-2}\phi_i\equiv 1 \,\,{\rm on}\,\,\Omega_i; \,\,0\leq\phi_i\leq 1 \,\,{\rm on}\,\, \overline\Sigma,\,\,\forall 1\leq i\leq k, \ee
 and for some positive constant $b_1$ depending only on $b_0$ and the metric $g$,
 \be\label{phi-3-2}|\nabla_g\phi_i|\leq b_1,\quad\forall 1\leq i\leq k.\ee
 For any $u\in W^{1,2}(\Sigma)$ satisfying \eqref{comp-2}, we have $\phi_i u\in W^{1,2}(\Omega)$ for all $1\leq i\leq k$, and thus
 \eqref{T-M-weak-2} implies
  \bna
  \int_\Sigma e^udv_g&\leq&\f{1}{\gamma_0}\int_{\Omega_i}e^{u}dv_g\\
  &\leq&\f{1}{\gamma_0}\int_{\Sigma}e^{\phi_iu}dv_g\\
  &\leq&\f{1}{\gamma_0}\exp\le(\f{1}{8\pi}\|\nabla_g(\phi_iu)\|_{L^2(\Sigma)}^2+\f{1}{|\Sigma|}\int_\Sigma \phi_iu dv_g +C\ri).
  \ena
  Note that \eqref{phi-1-2} gives
  $$\sum_{i=1}^k\|\nabla_g(\phi_i u)\|_{L^2(\Sigma)}^2=\le\|\nabla_g\le(u\sum_{i=1}^k\phi_i\ri)\ri\|_{L^2(\Sigma)}^2,$$
  and \eqref{phi-2-2} implies
  $$\sum_{i=1}^k\int_\Sigma \phi_iudv_g\leq\int_\Sigma |u|dv_g\leq \f{1}{|\Sigma|^{1/2}}\|u\|_{L^2(\Sigma)}. $$
  Combining the above three estimates,  \eqref{phi-3-2}, the Young inequality and an elementary inequality
  $$a\leq (a_1\cdots a_k)^{1/k}\,\, {\rm if}\,\,
  0\leq a\leq a_i, \, i=1,\cdots,k,$$ we obtain
  \bea\nonumber
  \int_\Sigma e^udv_g&\leq&\f{1}{\gamma_0}\le(\prod_{i=1}^k\exp\le(\f{1}{8\pi}\|\nabla_g(\phi_iu)\|_{L^2(\Sigma)}^2
  +\f{1}{|\Sigma|}\int_\Sigma \phi_iu dv_g+C\ri)\ri)^{1/k}\\
  \nonumber&=&\f{e^C}{\gamma_0}\exp\le(\f{1}{8k\pi}\sum_{i=1}^k\|\nabla_g(\phi_iu)\|_{L^2(\Sigma)}^2+
  \f{1}{k}{\f{1}{|\Sigma|}}\sum_{i=1}^k\int_\Sigma\phi_i udv_g\ri)\\\nonumber
  &=&\f{e^C}{\gamma_0}\exp\le(\f{1}{8k\pi}\le\|\nabla_g\le(u\sum_{i=1}^k\phi_i\ri)\ri\|_{L^2(\Sigma)}^2+
  \f{1}{k}{\f{1}{|\Sigma|}}\sum_{i=1}^k\int_\Sigma\phi_i udv_g\ri)\\\label{uL2-2}
  &\leq& C\exp\le(\f{1}{8k\pi}(1+\epsilon_1)\|\nabla_gu\|_{L^2(\Sigma)}^2+C(\epsilon_1)\|u\|_{L^2(\Sigma)}^2\ri).
  \eea

  Let $0<\lambda_1\leq \lambda_2\leq \cdots\leq \lambda_\ell\leq\lambda_{\ell+1}\leq\cdots$ be all eigenvalues of the
  Laplace-Beltrami operator with respect to the Neumann boundary condition. Clearly $\lambda_i\ra+\infty$ as $i\ra\infty$. Let $\{e_i\}_{i=1}^\infty$ be the corresponding
  unit normal eigenfunctions, i.e., $\Delta_g e_i=\lambda_ie_i$, $\int_\Sigma e_idv_g=0$, $\int_\Sigma e_ie_jdv_g=\delta_{ij}$ for $i,j=1,2,\cdots$. It is known that $\mathscr{H}:=\{u\in W^{1,2}(\Sigma): \overline{u}=0\}=E_\ell\oplus E_\ell^\perp$, where $E_\ell={\rm span}\{e_1,\cdots,e_\ell\}$ and $E_\ell^\perp={\rm span}\{e_{\ell+1},e_{\ell+2},
  \cdots\}$.
  Let $u\in W^{1,2}(\Sigma)$ be given as above. We decompose $u-\overline{u}=v+w$ with $v\in E_\ell$ and $w\in E_\ell^\perp$. Thus the Poincar\'e inequality implies
  $$\|v\|_{C^0(\overline\Sigma)}\leq \sum_{i=1}^\ell \|e_i\|_{C^0(\overline\Sigma)}\int_\Sigma |u-\overline{u}||e_i|dv_g\leq C_\ell\|\nabla_gu\|_{L^2(\Sigma)}.$$
  While by the definition of the $(\ell+1)$-th eigenvalue $\lambda_{\ell+1}$,
  $$\int_\Sigma w^2dv_g\leq \f{1}{\lambda_{\ell+1}}\int_\Sigma|\nabla_gw|^2dv_g.$$
  Having the above two estimates and applying \eqref{uL2-2} to $w$, we have
  \bna
  \int_\Sigma e^{u-\overline{u}}dv_g&\leq& e^{C_\ell\|\nabla_gu\|_{L^2(\Sigma)}}\int_\Sigma e^wdv_g\\
  &\leq& Ce^{C_\ell\|\nabla_gu\|_{L^2(\Sigma)}}\exp\le(\f{1}{8k\pi}(1+\epsilon_1)\|\nabla_gw\|_{L^2(\Sigma)}^2+
  \f{C(\epsilon_1)}{\lambda_{\ell+1}}\|\nabla_gw\|_{L^2(\Sigma)}^2\ri)\\
  &\leq&Ce^{C_\ell\|\nabla_gu\|_{L^2(\Sigma)}}\exp\le(\f{1}{8k\pi}\le(1+\epsilon_1+\f{C(\epsilon_1)}{\lambda_{\ell+1}}\ri)
  \|\nabla_gu\|_{L^2(\Sigma)}^2\ri).
  \ena
  This together with the Young inequality gives
  \be\label{sec-2}\log\int_\Sigma e^{u-\overline{u}}dv_g\leq \f{1}{8k\pi}\le(1+\epsilon_1\f{C(\epsilon_1)}{\lambda_{\ell+1}}+\epsilon_1\ri)
  \int_\Sigma |\nabla_gu|^2dv_g+C_{\ell,k,\epsilon_1}.
 \ee
 Let $0<\epsilon<8k\pi$ be any given number. Choosing $\epsilon_1=\epsilon/(16k\pi-2\epsilon)$, and then taking a sufficiently large $\ell$ such that  ${C(\epsilon_1)}/{\lambda_{\ell+1}}\leq 1$, we have by \eqref{sec-2} that
 $$\log\int_\Sigma e^{u-\overline{u}}dv_g\leq \f{1}{8k\pi-\epsilon}\int_\Sigma |\nabla_gu|^2dv_g+C,$$
 where $C$ is a constant depending only on $b_0$, $\gamma_0$, $k$, and $\epsilon$. This is exactly \eqref{imp-2}. $\hfill\Box$\\

 Define a functional $J_\rho:W^{1,2}(\Sigma)\ra\mathbb{R}$ by
 \be\label{funct-2}J_\rho(u)=\f{1}{2}\int_\Sigma|\nabla_gu|^2dv_g-\rho\log\int_\Sigma he^{u}dv_g+\f{\rho}{|\Sigma|}\int_\Sigma udv_g.\ee
 From now on to the end of this section, $J_\rho$ is always given as in \eqref{funct-2}.

\subsection{Continuous maps between sub-levels of \texorpdfstring{$J_\rho$}{the functional} and \texorpdfstring{$\Sigma_k$}{baycenter}}\label{3.3}

Let $\overline\Sigma_k$ be defined as in  \eqref{O-S-k}, and $(\p\Sigma)_k$ be defined by
$$(\p\Sigma)_k=\le\{\sum_{i=1}^kt_i\delta_{x_i}:t_i\geq 0,\sum_{i=1}^k=1,x_i\in\p\Sigma\ri\}.$$

\begin{lemma}\label{Psi-2}
Let $\rho\in(4k\pi,4(k+1)\pi)$. Then for any sufficiently large $L>0$, there exists a continuous retraction
$$\Psi:J_\rho^{-L}=\le\{u\in W^{1,2}(\Sigma): J_\rho(u)\leq -L\ri\}\ra \overline{\Sigma}_k.$$
Moreover, if $(u_n)\subset W^{1,2}(\Sigma)$ satisfies
$\f{e^{u_n}}{\int_\Sigma e^{u_n}dv_g}dv_g\ra\sigma\in\overline\Sigma_k$, then $\Psi(u_n)\ra \sigma\in\overline\Sigma_k$.
\end{lemma}
\proof Since the proof is almost the same as that of Lemma \ref{Psi}, we omit the details here.  $\hfill\Box$\\

For any finite set $E$, we denote the number of all distinct points of $E$ by $\sharp E$. We define
\be\label{subset}
\mathscr{S}_k=\left\{\sigma\in\overline{\Sigma}_k: \sharp{({\rm supp}\sigma\cap\Sigma)}+
\sharp{{\rm supp}\sigma}\leq k\right\}.
\ee
Let us explain what $\mathscr{S}_k$ means. Clearly, $\mathscr{S}_1=(\p\Sigma)_1=\p\Sigma$; For $k=2$, $\mathscr{S}_2=\{\delta_x:x\in\Sigma\}\cup\{t\delta_{x_1}+(1-t)\delta_{x_2}:0\leq t\leq 1,x_1,x_2\in\p\Sigma\}=(\p\Sigma)_2\cup\Sigma$; For $k=3$, we write
$\mathscr{S}_3=\{t\delta_{x_1}+(1-t)\delta_{x_2}:0\leq t\leq 1, x_1\in\Sigma,x_2\in\p\Sigma\}\cup\{t_1\delta_{x_1}+t_2\delta_{x_2}+t_3
\delta_{x_3}: 0\leq t_i\leq 1, x_i\in\p\Sigma,1\leq i\leq 3, t_1+t_2+t_3=1\}=(\p\Sigma)_3\cup \mathscr{A}_3$, where
$\mathscr{A}_3=\{t\delta_{x_1}+(1-t)\delta_{x_2}:0\leq t\leq 1, x_1\in\Sigma,x_2\in\p\Sigma\}$. We observe that ${\rm dim}\mathscr{A}_3<
{\rm dim}(\p\Sigma)_3=5$, since $(\p\Sigma)_3\setminus (\p\Sigma)_2$ is a smooth $5$-dimensional manifold, and
${\rm dim}\mathscr{A}_3\leq 4$.
\begin{lemma}\label{Sk-nonc}
For any $k\geq 1$, $\mathscr{S}_k$ is non-contractible.
\end{lemma}
\proof Obviously $\mathscr{S}_1=\p\Sigma$ is non-contractible. For $k\geq 2$, based on the above observation, an induction argument shows
$\mathscr{S}_k=(\p\Sigma)_k\cup\mathscr{A}_k$,
where ${\rm dim}\mathscr{A}_k<{\rm dim}(\p\Sigma)_k=2k-1$. Though $(\p\Sigma)_k$ is a combination of several
different dimensional branches, we still denote the maximum dimension of those branches by
${\rm dim}(\p\Sigma)_k$. Arguing as \cite[Lemma 4.7]{Djadli}, we have that $(\p\Sigma)_k$ is non-contractible. (This was also noticed by Zhang-Zhou-Zhou \cite{Zhang-Zhou-Zhou}). In fact, we have
 \begin{equation}\label{non-00}
  H_{2k-1}((\p\Sigma)_k,\mathbb{Z}_2)\not=\{0\}.  
 \end{equation}
 Given any $(2k-1)$-dimensional closed chain $\mathcal{C}_{2k-1}$ and any $(2k-1)$-dimensional boundary chain ${\mathcal{E}}_{2k-1}$ of $(\p\Sigma)_k$.
 Since $(\p\Sigma)_k$ is a closed sub-topological space in $\mathscr{S}_k$, and ${\rm dim}((\p\Sigma)_k)=
 {\rm dim}(\mathscr{S}_k)=2k-1$, one easily sees that $\mathcal{C}_{2k-1}$ is also a closed chain of $\mathscr{S}_k$
 and $\mathcal{E}_{2k-1}=0$ is also the boundary $(2k-1)$-chain of $\mathscr{S}_k$. Hence  $$H_{2k-1}((\p\Sigma)_k,\mathbb{Z}_2)\subset H_{2k-1}(\mathscr{S}_k,\mathbb{Z}_2),$$ which together with (\ref{non-00}) implies
 that
$\mathscr{S}_k$ is non-contractible. $\hfill\Box$  \\

 Take a smooth increasing function $\eta:\mathbb{R}\ra\mathbb{R}$
satisfying  $\eta(t)=t$ for $t\leq 1$, and $\eta(t)=2$ for $t\geq 2$. Set $\eta_r(t)=r\eta(t/r)$ for $r>0$.
For $\lambda>0$, $x\in\overline\Sigma$, and $1\leq\ell<m$, we define
\be\label{ppp-2}\widetilde{\phi}_{\lambda,\sigma}(x)=\log\le(\sum_{i=1}^\ell \f{t_i}{2}\f{8\lambda^2}{(1+\lambda^2\eta_r^2
({\rm dist}(x,x_i)))^2}+\sum_{i=\ell+1}^m {t_i}\f{8\lambda^2}{(1+\lambda^2\eta_r^2
({\rm dist}(x,x_i)))^2}\ri)\ee
and
\be\label{p-1-2}\phi_{\lambda,\sigma}(x)=\widetilde{\phi}_{\lambda,\sigma}(x)-\f{1}{|\Sigma|}\int_\Sigma\widetilde{\phi}_{\lambda,\sigma}dv_g.\ee

\begin{lemma}\label{P-2}
Let $\rho\in(4k\pi,4(k+1)\pi)$. If $\lambda>0$ is chosen sufficiently large, and $r>0$ is chosen sufficiently small, then for any $\sigma\in\mathscr{S}_k$, there hold
\be\label{J-less-2}J_\rho(\phi_{\lambda,\sigma})\leq \le(4k\pi-\rho\ri)\log\lambda\ee
and
\be\label{tend-2}\f{e^{\phi_{\lambda,\sigma}}}{\int_\Sigma e^{\phi_{\lambda,\sigma}}dv_g}dv_g\ra \sigma\quad{\rm as}\quad\lambda\ra+\infty.\ee
\end{lemma}
\proof Both the cases ${\rm supp}\sigma\cap\Sigma=\varnothing$ and ${\rm supp}\sigma\cap\Sigma\not=\varnothing$ can be dealt with in the same way. Given $\sigma\in \mathscr{S}_k$. With no loss of generality, we assume ${\rm supp}\sigma=\{x_1,\cdots,x_m\}\subset\overline\Sigma$,
${\rm supp}\sigma\cap \Sigma=\{x_1,\cdots,x_\ell\}$, and $m+\ell\leq k$.
 Let $\widetilde{\phi}_{\lambda,\sigma}$ and $\phi_{\lambda,\sigma}$ be defined as in \eqref{ppp-2}
and \eqref{p-1-2} respectively, where $\lambda>0$ and $r>0$.
Write $r_i=r_i(x)={\rm dist}(x,x_i)$ for $x\in\overline\Sigma$. A simple observation gives
\be\label{wide-2}
\widetilde{\phi}_{\lambda, \sigma}(x)=\left\{\begin{array}{lll}
\log \frac{8 \lambda^2}{(1+4 \lambda^{2}r^{2} )^2 } &{\rm for}& x \in \overline\Sigma\setminus \cup_{i=1}^{k} B_{2 r}(x_i)
\\[1.5ex]
\log \left(\frac{8 \lambda^2 \overline{t}_{i}}{ \left(1+\lambda^{2} \eta_r^{2}\left(r_i\right) \right)^2 }+\frac{8 \lambda^2\left(1-\overline{t}_{i}\right)}{ (1+4\lambda^{2} r^{2})^2}\right) &{\rm for}& x \in B_{2 r}(x_i),
\end{array}\right.
\ee
where $\overline{t}_i=t_i/2$ for $1\leq i\leq\ell$, $\overline{t}_i=t_i$ for $\ell+1\leq i\leq m$,  and $B_{2 r}(x_i)=\{x\in\overline\Sigma:{\rm dist}(x,x_i)<2r\}$ denotes a geodesic ball centered at $x_i$ with radius $2r$.
One easily sees  $\phi_{\lambda,\sigma}\in W^{1,2}(\Sigma)$ and $\int_\Sigma\phi_{\lambda,\sigma}dv_g=0$.
 For $x\in B_{2r}(x_i)$, $i=1,\cdots,m$, a straightforward calculation shows
$$\nabla_g\widetilde\phi_{\lambda,\sigma}(x)=\f{\f{8\lambda^2\overline{t}_i}{(1+\lambda^2\eta_r^2(r_i))^2}}
{\f{8\lambda^2\overline{t}_i}{(1+\lambda^2\eta_r^2(r_i))^2}+\f{8\lambda^2(1-\overline{t}_i)}{(1+4\lambda^2r^2)^2}}
\f{4\lambda^2\eta_r(r_i)\eta_r^\prime(r_i)\nabla_gr_i}{1+\lambda^2\eta_r^2(r_i)},$$
and thus
$$|\nabla_g\widetilde{\phi}_{\lambda,\sigma}(x)|\leq \f{4\lambda^2\eta_r(r_i)\eta_r^\prime(r_i)}{1+\lambda^2\eta_r^2(r_i)}.$$
In view of \eqref{wide-2}, there holds $\nabla_g\widetilde{\phi}_{\lambda,\sigma}(x)=0$ for $x\in \overline\Sigma\setminus \cup_{i=1}^{m} B_{2 r}(x_i)$.
We calculate for $1\leq i\leq \ell$,
$$\int_{B_{2r(x_i)}}
\le(\f{4\lambda^2\eta_r(r_i)\eta_r^\prime(r_i)}{1+\lambda^2\eta_r^2(r_i)}\ri)^2dv_g=16\pi(1+O(r^2))\le(\log(1+\lambda^2r^2)
+\f{1}{1+\lambda^2r^2}-1\ri);$$
and for $\ell+1\leq i\leq m$,
$$\int_{B_{2r(x_i)}}
\le(\f{4\lambda^2\eta_r(r_i)\eta_r^\prime(r_i)}{1+\lambda^2\eta_r^2(r_i)}\ri)^2dv_g=8\pi(1+O(r^2))\le(\log(1+\lambda^2r^2)
+\f{1}{1+\lambda^2r^2}-1\ri).$$
For a fixed $r>0$, since $x_1,\cdots,x_m$ are arbitrary, one sees that
$\{B_{2r}(x_1),\cdots,B_{2r}(x_m)\}$ may have nonempty intersections,  and thus 
\bea
\nonumber \int_\Sigma |\nabla_g\widetilde{\phi}_{\lambda,\sigma}|^2dv_g&=&\int_{\cup_{i=1}^mB_{2r}(x_i)}|\nabla_g\widetilde{\phi}_{\lambda,\sigma}|^2dv_g\\
\nonumber&\leq&\sum_{i=1}^m\int_{B_{2r(x_i)}}
\le(\f{4\lambda^2\eta_r(r_i)\eta_r^\prime(r_i)}{1+\lambda^2\eta_r^2(r_i)}\ri)^2dv_g\\
\nonumber&\leq&\le((16\pi\ell+8\pi(m-\ell)\ri)(1+O(r^2))\le(\log(1+\lambda^2r^2)+\f{1}{1+\lambda^2r^2}-1\ri)+O(1)\\
&\leq& 8k\pi(1+O(r^2))\log \lambda^2+C\label{nabl-2}
\eea
for some constant $C$ depending only on $r$. Moreover, for any $s$, $0<s<\min\{r,\f{1}{2}\min_{1\leq i< j\leq m}{\rm dist}(x_i,x_j)\}$, there holds
\bna
\int_{\cup_{i=1}^mB_{2r}(x_i)}e^{\widetilde{\phi}_{\lambda,\sigma}}dv_g&=&\int_{\cup_{i=1}^mB_{s}(x_i)}e^{\widetilde{\phi}_{\lambda,\sigma}}dv_g
+O\le(\f{1}{\lambda^2s^2}\ri)\\
&=&\sum_{i=1}^m\int_{B_{s}(x_i)}
\f{8\lambda^2\overline{t}_i}{(1+\lambda^2r_i^2)^2}dv_g+O\le(\f{1}{\lambda^2s^2}\ri)\\
&=&4\pi(1+O(s^2))+O\le(\f{1}{\lambda^2s^2}\ri),
\ena
and
\bna
\int_{\Sigma\setminus\cup_{i=1}^mB_{2r}(x_i)}e^{\widetilde{\phi}_{\lambda,\sigma}}dv_g=O\le(\f{1}{\lambda^2r^4}\ri).
\ena
It follows that
\be\label{sense-2}\int_\Sigma e^{\widetilde\phi_{\lambda,\sigma}}dv_g=4\pi(1+O(s^2))+O\le(\f{1}{\lambda^2s^2}\ri)+O\le(\f{1}{\lambda^2r^4}\ri).\ee
Passing to the limit $\lambda\ra+\infty$ first, and then $s\ra 0+$, we have
\be\label{limit-2}\lim_{\lambda\ra+\infty}\int_\Sigma e^{\widetilde\phi_{\lambda,\sigma}}dv_g=4\pi.\ee
Note that there exists some constant $C$ depending only on $r$ such that
\be\label{ggg-2}\f{1}{|\Sigma|}\int_\Sigma \widetilde\phi_{\lambda,\sigma}dv_g\leq -\log\lambda^2+C.\ee
Hence by \eqref{limit-2} and \eqref{ggg-2},
\bea\nonumber\int_\Sigma e^{\phi_{\lambda,\sigma}}dv_g&\geq& C\le(1+o_\lambda(1)\ri)\lambda^2.
\eea
This together with \eqref{nabl-2} gives
\bna
J_\rho(\phi_{\lambda,\sigma})&=&\f{1}{2}\int_\Sigma|\nabla_g\phi_{\lambda,\sigma}|^2dv_g-\rho\log\int_\Sigma he^{\phi_{\lambda,\sigma}}dv_g\\
&\leq&(4k\pi-\rho+O(r^2))\log\lambda^2+C_r.
\ena
Since $\rho>4k\pi$, choosing $r>0$ sufficiently small and $\lambda>0$ sufficiently large, we conclude \eqref{J-less-2}.

Finally we prove \eqref{tend-2}. Let $\sigma=\sum_{i=1}^mt_i\delta_{x_i}\in\mathscr{S}_k$ be as above. For any $\varphi\in C^1(\overline\Sigma)$, similar to \eqref{sense-2}, we calculate
$$
\int_\Sigma\varphi e^{\widetilde{\phi}_{\lambda,\sigma}}dv_g
=4\pi\sum_{i=1}^kt_i\varphi(x_i)+O(s^2)+O\le(\f{1}{\lambda^2s^2}\ri)+O\le(\f{1}{\lambda^2r^4}\ri).
$$
Letting $\lambda\ra+\infty$ first, and then $s\ra 0+$, we get
$$\lim_{\lambda\ra+\infty}\int_\Sigma\varphi e^{\widetilde{\phi}_{\lambda,\sigma}}dv_g=4\pi\sum_{i=1}^kt_i\varphi(x_i).$$
This together with \eqref{limit-2} implies \eqref{tend-2}.
 $\hfill\Box$\\

 Similar to \eqref{proj}, for a sufficiently small $\epsilon_0>0$, we have a continuous retraction
 $$\mathfrak{p}:\{\sigma\in\mathcal{D}(\Sigma):{\mathbf d}(\sigma,\mathscr{S}_k)<\epsilon_0\}\ra\mathscr{S}_k.$$

 \begin{lemma}\label{PHI-2}
 Let $\Psi$ and $L$ be as in Lemma \ref{Psi-2}. If $\lambda>0$ is chosen sufficiently large,
 then there exists a continuous map $\Phi_\lambda: \mathscr{S}_{k}\ra J_\rho^{-L}$ such that
 $\mathfrak{p}\circ\Psi\circ \Phi_\lambda:\mathscr{S}_k\ra\mathscr{S}_k$ is homotopic to the identity map ${\rm Id}:\mathscr{S}_k\ra\mathscr{S}_k$.
 \end{lemma}

 \proof Let $\phi_{\lambda,\sigma}$ be constructed as in Lemma \ref{P-2}. For any $\sigma\in\mathscr{S}_k$, we define
 $\Phi_\lambda(\sigma)=\phi_{\lambda,\sigma}$ for large $\lambda>0$. Clearly the map $\Phi_\lambda:\mathscr{S}_k\ra W^{1,2}(\Sigma)$ is continuous.
 By \eqref{J-less-2}, if $\lambda\geq e^{L/(\rho-4k\pi)}$, then $J_\rho(\phi_{\lambda,\sigma})\leq -L$. Thus
 $\Phi_\lambda(\sigma)\in J_\rho^{-L}$. By Lemma \ref{Psi-2} and \eqref{tend-2}, there holds
 \bna
 \mathfrak{p}\circ\Psi\circ\Phi_\lambda(\sigma)&=&\mathfrak{p}\circ\Psi(\phi_{\lambda,\sigma})\\
 &=&\mathfrak{p}\circ\psi_k\le(\f{e^{\phi_{\lambda,\sigma}}}{\int_\Sigma e^{\phi_{\lambda,\sigma}}dv_g}dv_g\ri)\\
 &\ra&\sigma
 \ena
 as $\lambda\ra+\infty$. Hence $\mathfrak{p}\circ\Psi\circ \Phi_\lambda$ is homotopic to ${\rm Id}:\mathscr{S}_k\ra\mathscr{S}_k$. $\hfill\Box$

\subsection{Min-max values}\label{3.4}
 Let $$\mathcal{H}=\le\{u\in W^{1,2}(\Sigma):\int_\Sigma udv_g=0\ri\}$$
 and
$$\widehat{\mathscr{S}}_k=\mathscr{S}_k\times [0,1]/(\mathscr{S}_k\times\{0\})$$ be the topological cone over
$\mathscr{S}_k$. A path set associated to
the metric space
$\widehat{\mathscr{S}}_k$ is defined by
$$\Gamma_\lambda=\le\{\gamma\in C^0(\widehat{\mathscr{S}}_k,\mathcal{H}):\gamma|_{\widehat{\mathscr{S}}_k\times\{1\}}\in\Gamma_{\lambda,0}\ri\},$$
where $\Gamma_{\lambda,0}$ is given by
$$\Gamma_{\lambda,0}=\le\{\gamma\in C^0({\mathscr{S}}_k\times\{1\},\mathcal{H}):\gamma(\sigma,1)=\Phi_\lambda(\sigma),\forall \sigma\in{\mathscr{S}}_k\,\ri\}.$$
If we write a path $\overline{\gamma}:\widehat{\mathscr{S}}_k\ra \mathcal{H}$ by $\overline{\gamma}(\sigma,t)=t\phi_{\lambda,\sigma}$,
then $\overline{\gamma}\in \Gamma_\lambda$, and thus $\Gamma_\lambda\not=\varnothing$.

For real numbers $\lambda$ and $\rho$, we set
$$\alpha_{\lambda,\rho}=\inf_{\gamma\in\Gamma_\lambda}\sup_{(\sigma,t)\in\widehat{\mathscr{S}}_k}J_\rho(\gamma(\sigma,t))$$
and
$$\beta_{\lambda,\rho}=\sup_{\gamma\in\Gamma_{\lambda,0}}\sup_{(\sigma,t)\in{\mathscr{S}}_k\times\{1\}}
J_\rho(\gamma(\sigma,t)).$$
\begin{lemma}
Let $\rho\in(4k\pi,4(k+1)\pi)$.
If $\lambda$ is chosen sufficiently large, and $r$ is chosen sufficiently small, then
$-\infty<\beta_{\lambda,\rho}<\alpha_{\lambda,\rho}<+\infty$.
\end{lemma}
\proof The proof is very similar to that of Lemma \ref{mini}. It suffices to use Lemma \ref{PHI-2} instead of
the fact $\pi\circ\Psi\circ\Phi_\lambda:\Sigma_{\epsilon,k}\ra\Sigma_{\epsilon,k}$ is homotopic to ${\rm Id}:
\Sigma_{\epsilon,k}\ra\Sigma_{\epsilon,k}$, and use Lemma \ref{Sk-nonc} instead of the non-contractibility of $\Sigma_{\epsilon,k}$.
 $\hfill\Box$

 \subsection{Completion of the proof of \autoref{thm2}}\label{3.5}
 Denote $\alpha_\rho=\alpha_{\lambda,\rho}$ for sufficiently large $\lambda>0$.
 Similar to Lemma \ref{monotonicity}, $\alpha_\rho/\rho$ is decreasing in $\rho\in (4k\pi,4(k+1)\pi)$.
 Let $$\Lambda_{k}=\le\{\rho\in(4k\pi,4(k+1)\pi): \f{\alpha_\rho}{\rho}\,\,{\rm is\,\, differentiable\,\, at}\,\,\rho\ri\}.$$
 In view of an analog of Lemma \ref{dense}, $\alpha_\rho$ is a critical value of $J_\rho$ for any $\rho\in \Lambda_k$.

 Now we let $\rho\in(4k\pi,4(k+1)\pi)$.
    Take an increasing sequence of numbers
  $(\rho_n)\subset\Lambda_{k}$ such that  $\rho_n\ra \rho$, $(\rho_n)\subset[a,b]\subset(4k\pi,4(k+1)\pi)$,
  and $\alpha_{\rho_n}$ is achieved by $u_n\in \mathcal{H}$. Moreover,
  $u_n$ satisfies the Euler-Lagrange equation
  \be\label{n-equ-2}\Delta_gu_n=\rho_n\le(\f{he^{u_n}}{\int_\Sigma he^{u_n}dv_g}-\f{1}{|\Sigma|}\ri).\ee
  Since $\alpha_\rho/\rho$ is decreasing in $\rho\in [a,b]$,
  \be\label{alpha-n-2}\alpha_{\rho_n}\leq \f{b}{a}\alpha_a.\ee
  Denoting $v_n=u_n-\log\int_\Sigma he^{u_n}dv_g$, we have
  $$\le\{\begin{array}{lll}
  \Delta_gv_n=\rho_n(he^{v_n}-|\Sigma|^{-1})\\[1.5ex]
  \int_\Sigma he^{v_n}dv_g=1.
  \end{array}\ri.$$
  By Lemma \ref{compact-2}, $(u_n)$ is bounded in $L^\infty(\overline{\Sigma})$. Let $\Omega_1,\cdots,\Omega_{k+1}$ be disjoint
  closed sub-domains
  of $\overline\Sigma$. It follows from Lemma \ref{improved-2} that
  $$\log\int_\Sigma e^{u_n}dv_g\leq \f{1}{8(k+1)\pi-\epsilon}\int_\Sigma|\nabla_gu_n|^2dv_g+C_\epsilon$$
  for any $\epsilon>0$ and some constant $C_\epsilon>0$. This together with \eqref{alpha-n-2} implies that for $0<\epsilon<8(k+1)\pi-2b$,
  \bna
  \f{1}{2}\int_\Sigma|\nabla_gu_n|^2dv_g&=&J_{\rho_n}(u_n)+\rho_n\log\int_\Sigma he^{u_n}dv_g\\
  &\leq&\f{b}{8(k+1)\pi-\epsilon}\int_\Sigma|\nabla_gu_n|^2dv_g+C.
  \ena
  Then it follows that $(u_n)$ is bounded in $\mathcal{H}$. With no loss of generality, we assume
  $u_n$ converges to $u_0$ weakly in $\mathcal{H}$, strongly in $L^p(\Sigma)$ for any $p>1$, and almost everywhere in $\overline\Sigma$.
  Moreover, $e^{u_n}$ converges to $e^{u_0}$ strongly in $L^p(\Sigma)$ for  any $p>1$. By \eqref{n-equ-2}, $u_0$ satisfies
    $$\Delta_gu_0=\rho\le(\f{he^{u_0}}{\int_\Sigma he^{u_0}dv_g}-\f{1}{|\Sigma|}\ri)$$
    in the distributional sense.
  In particular $u_0$ is a critical point of $J_\rho$. $\hfill\Box$

\hspace{3cm}\\


\end{document}